\documentclass[11pt,leqeq]{article}
\usepackage{amssymb,amsthm}
\usepackage{amsmath}
\usepackage{amscd}
\usepackage{xy}
\xyoption{all}
\usepackage[dvips]{graphicx}
\newtheorem{thm}{Theorem}%[section] (If you want theorem numbered
\newtheorem{lema}[thm]{Lemma}%               with section number.  Same
\newtheorem{cor}[thm]{Corollary}%       goes for lemmas, etc.)
\newtheorem{prop}[thm]{Proposition} %--> \begin\end{theorem,lemma,...}

\newtheorem{defi}[thm]{Definition}

\DeclareFontFamily{OMS}{rsfs}{\skewchar\font'60}
\DeclareFontShape{OMS}{rsfs}{m}{n}{<-5>rsfs5 <5-7>rsfs7 <7->rsfs10
}{} \DeclareSymbolFont{rsfs}{OMS}{rsfs}{m}{n}
\DeclareSymbolFontAlphabet{\scr}{rsfs}

\newcommand{\des}{\displaystyle}
\newcommand{\im}{{\rm{im}}}

\newcommand{\ma}{{\rm Map}}
\newcommand{\Co}{{\rm Cob}}

\newcommand{\Mfu}{{\rm MFunc}}

\newcommand{\Ho}{{\rm H}}

\newcommand{\Ca}{{\rm C}}

\newcommand{\ri}{{\rm r}}
\newcommand{\ve}{{\rm vect}}

\newcommand{\de}{{\rm dg\mbox{-}vect}}

\newcommand{\ob}{{\rm Ob}}
\newcommand{\HL}{{\rm HLQFT}}

\newcommand{\id}{{\rm id}}
\newcommand{\ho}{{\rm Hom}}
\newcommand{\sC}{\scr{C}}
\newcommand{\OB}{{\rm ob}}
\newcommand{\I}{{\rm I}}
\renewcommand{\des}{\displaystyle}
\def\N{{\Bbb N}}

\def\O{{\cal{O}}}

\def\C{{\cal{C}}}

\def\d{{\cal{D}}}
\def\f{{\cal{F}}}

\def\h{{\cal{H}}}

\date{}

\begin{document} % End of preamble and beginning of text.

\title{Homological Quantum Field Theory} % Declares the document's title.
\author{Edmundo  Castillo and Rafael D\'\i az }
\maketitle

\begin{abstract}
We show that the space of chains of smooth maps from spheres into
a fixed compact oriented manifold has a natural structure of a
transversal $d$-algebra.  We construct a structure of transversal
$1$-category on the space of chains of maps from a suspension
space $S(Y)$, with certain boundary restrictions, into a fixed
compact oriented manifold. We define homological quantum field
theories \HL \ and construct several examples of such structures.
Our definition is based on the notions of string topology of Chas
and Sullivan, and homotopy quantum field theories of Turaev.
\end{abstract}

\section{Introduction}

This work takes part in the efforts aimed to understand the
mathematical foundations of quantum field theory by unveiling its
underlying categorical structures. A distinctive feature of the
categorical approach to field theory, to be reviewed in Section
\ref{sec7}, is that it works better for theories with a rather
large group of symmetries, for example, a theory invariant under
arbitrary topological transformations. Thus our subject matter is
deeply intertwined with algebraic topology. A major problem in
algebraic topology with a long and rich history
\cite{FM, Mil, GP,  SSmale, WPT} is the classification of compact smooth manifolds, i.e.
the description of the set of equivalences classes of manifolds,
where two manifolds are equivalent if they are diffeomorphic. It
has been proven that a classification is possible in all
dimensions in principle except for the case dimension $4$ which
remains open. The $3$ dimensional case have been settled  with the
completion by Pereleman of the Ricci flow approach to the Thurston
geometrization program. However only in dimension two we can
distinguish manifolds in an efficient computable way, i.e. given a
couple of manifolds of the same dimension it is a hard problem to
tell whether they are diffeomorphic manifolds or not. The standard
way to distinguish non-diffeomorphic manifolds $M_1$ and $M_2$ is
to find a topological invariant $I$ such that $I(M_1)\neq I(M_2).$
By definition a topological invariant with values in a ring $R$ is
a map that assigns an element $I(M)$ of $R$ to each manifold $M$
in such a way that diffeomorphic manifolds are mapped into the
same element. Topological invariants thus provide a way to
effectively distinguish non-diffeomorphic manifolds. A
quintessential example is the Euler characteristic $\chi$, a
topological  invariant powerful enough to classify $2$ dimensional
manifolds. An early achievement in algebraic topology, one whose
consequences via the development of the theory of categories of
Eilenberg and Mac Lane \cite{SMacLane, SMacLane2} has already
reshaped modern mathematics, was the realization of the  necessity
to study not only ring valued invariants, but also invariants
taking values in arbitrary categories. Thus in this more general
approach a topological invariant is a functor from the category of
compact manifolds into a fixed target category. A prominent
example is singular homology $H$, the functor that assigns  to
each manifold $M$ the homology $H(M)=H(C(M))$  of the
$\mathbb{N}$-graded vector space $C(M)$ of singular chains on $M$,
i.e. the space of closed singular chains modulo exact ones. For
simplicity we always consider homology with complex coefficients.
The Euler characteristic is the super-dimension of $H(M)$, i.e. it
is given by
$$\chi(M)=\mathrm{sdim}(H(M))=\sum_{i=0}^{\dim(M)}(-1)^i\dim(H^i(M)).$$
Hence homology is a categorification of the Euler characteristic,
i.e. it is a functor with values in graded vector spaces linked to
the Euler characteristic via the notion of super-dimension. The
process of categorification is nowadays under active research from
a wide range  of viewpoints, see for example
\cite{PP, BaezDolan, Blan1, Blan2, Blan3, Cas8, cas2, cy, RDEP, RDEP1, kho}.
The example above illustrates the essence of the categorification
idea: there are manifolds with equal Euler characteristic but
non-isomorphic singular homology groups, thus the latter invariant
is subtler and deeper. Likewise, in contrast with the
correspondence
$$M \longrightarrow H(M)$$ which has been extensively studied in the
literature, the correspondence $$ M
\longrightarrow C(M),$$ is much less understood, although its core properties has
been elucidated by Mandell in his works \cite{mand, mand2}. The
latter correspondence contains deeper information for the same
reasons that homology contains deeper information than the Euler
characteristic: a further level of categorification has been
achieved, or in more mundane terms, non quasi-isomorphic complexes
may very well have isomorphic homology. We are going to take this
rather subtle issue seriously and make an effort to work
consistently at the chain level, rather than at the purely
homological level. This will require, among other things, that we
use a different model for homology in place of singular homology.
Chains on a manifold $M$ in
our model are smooth maps from manifolds with corners into $M$. \\

In this work we study topological invariants  for compact oriented
manifolds coming from the following simple idea. Fix a compact
manifold $L$ and for each manifold $M$ consider the topological
space
$$M^L= \{x \mid x \colon L \longrightarrow M \mbox{\ \ piecewise smooth map } \}$$
provided with the compact-open topology. The map sending $M$ into
the homology $H(M^L)$ of $M^L$ is a topological invariant which
assigns to each manifold a $\mathbb{N}$-graded vector space. An
interesting fact that will reemerge at various points in this work
is that if $L$ is chosen conveniently then the space $H(M^L)$ is
naturally endowed with a rich algebraic structure. Let us
highlight a few landmarks in the historical development of this
fruitful idea. The first example comes from classical algebraic
topology. Given a topological space $M$ with a marked point $p \in
M$ consider the space $M_p^{S^1}$ of loops in $M$ based at $p$,
i.e. the space
$$M_p^{S^1}=\{ x \mid
x : S^{1}\longrightarrow M \ \mbox{ piecewise smooth and } \
\gamma(1)=p  \}$$ provided with the compact-open topology.  The
space of based loops comes with a natural product given by
concatenation
$$M_p^{S^1}\times M_p^{S^1}
\longrightarrow M_p^{S^1}.$$ This product, introduced by introduced by
Pontryagin,  is associative up to homotopy. By  the
K$\ddot{\mbox{u}}$nneth formula and functoriality of homology the
Pontryagin product induces an associative product
$$H(M_p^{S^1})\otimes
H(M_p^{S^1})\longrightarrow H(M_p^{S^1})$$ on the homology groups
of  $M_p^{S^1}.$ Stasheff  in his celebrated works
\cite{S1, S2} introduced $A_{\infty}$-spaces and
$A_{\infty}$-algebras as tools for the study of spaces
homotopically equivalent to topological monoids. The primordial
example of an $A_{\infty}$-space is precisely $M_p^{S^1}$ the
space of based loops. Likewise singular chains $C(M_p^{S^1})$ on
$M_p^{S^1}$ are the quintessential example of an
$A_{\infty}$-algebra. The $A_{\infty}$-structure on $C(M_p^{S^1})$
induces an associative product on the homology groups
$H(M_p^{S^1})$ which agrees with the Pontryagin product. In this
work we do not deal explicitly with $A_{\infty}$-algebras or
$A_{\infty}$-categories, instead we shall use $1$-algebras and
$1$-categories. However, the reader should be aware that these
notions
are, respectively, equivalent.\\

A second flow of ideas came from string theory, a branch of high
energy physics that has been proposed by a distinguished group of
physicist -- references \cite{W1,WZ,BZ}  are not too far from the
spirit of this work -- as a unifying theory for all fundamental
forces of nature, including the standard model of nature and
general relativity. The primordial object of study in string
theory is the dynamics of a small loop moving inside a manifold
$M$, i.e. in string theory the configuration space  $M$ is the
infinite-dimensional space
$$M^{S^1}=\{x \mid x \colon S^{1}\longrightarrow M \mbox{ smooth } \}$$
of non-based loops in $M$, provided with the compact-open
topology. The analytical difficulties present in string theory
have prevented, to this day, a fully rigorous mathematical
description. Chas and Sullivan in their seminal work
\cite{SCh}  initiated the study of strings using classical algebraic topology.
The key observation made by them is that even though $M^{S^1}$
does not posses a product analogue to the Pontryagin product, the
homology $H(M^{S^1})$ of $M^{S^1}$ comes with a natural
associative product, which generalizes the Goldman bracket
\cite{gold1, gold2} on homotopy classes of curves embedded in a compact Riemann
surface. It is natural to wonder if that product arises from a
product defined at the chain level. We have hit an important
subtlety that will be a major theme of this work: the fact that
the product at the chain level is naturally defined only for
transversal chains and only if we use an appropriated definition
of chains. To work with algebras, and more generally categories,
with a product defined only for transversal tuples, we shall adopt
the theory of transversal or partial algebras of Kriz and May
\cite{KM}. To define the product at the chain level in Section
\ref{sec1} we present a model for homology using manifolds with
corners instead of simplices as the possible domain for chains.
This construction is motivated by the observation that the
transversal intersection of chains with simplicial domain is in a
natural way a (sum of) chain(s) having as domain a manifold with
corners. With this provisions then one can show that indeed the
Chas-Sullivan product comes from an associative up to homotopy
product defined for transversal chains on the space of non-based
loops, more precisely, we show that the space of chains
is a transversal $1$-algebra.\\

Since its introduction the full range of structures  taking part
of string topology has been study and generalized from various
viewpoints, out which we mention just a few without pretension of
being exhaustive; for comprehensive reviews of string topology in
its various approaches the reader may consult
\cite{CoVo, S2}. In addition to the string product
there is a string bracket $$\{ \ , \ \}:\Ho(M^{S^1})\otimes
\Ho(M^{S^1})\longrightarrow \Ho(M^{S^1})$$ and a delta operator
$$\Delta:\Ho(M^{S^1})\longrightarrow \Ho(M^{S^1}),$$ which are defined in such a way
that they together with the string product give $\Ho(M^{S^1})$ the
structure of a Batalin-Vilkovisky or BV algebra. The space of
functionals of fields, including ghost and anti-ghost, of a gauge
theory is naturally endowed with the structure of a BV algebra
\cite{bat1, bat2}. Cattaneo, Fr$\ddot{\mbox{o}}$hlich and Pedrini have shown in
\cite{cat} that the bracket of the BV structure on $\Ho(M^{S^1})$
corresponds with the bracket of the BV structure on the
functionals of the higher dimensional Chern-Simons action
\cite{ale} with gauge group $GL(n,
\mathbb{C})$. Another interesting approach to string topology is obtained via
Hoschild cohomology, indeed Cohen and Jones show in
\cite{CJ} that there is a ring isomorphism
$$\Ho(M^{S^1}) \longrightarrow H(C^*(M), C^*(M)),$$
where $C^*(M)$ is the co-chain algebra of a simply connected
manifold $M$, and $H(C^*(M), C^*(M))$ is the Hoschild cohomology
of $C^*(M)$ the algebra of co-chains in $M$. The ring structure on
$H(C^*(M), C^*(M))$  is given by the Gerstenhaber cup product. It
turns out that this isomorphism preserves the full BV structure on
both sides, as shown in the recent works
\cite{YJM2, menich},  both based on the F\'elix, Thomas and
Vigu\'e-Porrier \cite{YJM1}  cochain model for the product on  $\Ho(M^{S^{1}})$ using tools from rational homotopy theory. \\

Let us mention three additional approaches to string topology.
Chataur in \cite{Ch} described string topology in terms of the
geometric cycles approach to homology \cite{J}. Cohen in
\cite{Ch2} studies string topology from the viewpoint of Morse
theory, and shows that the Floer homology $HF(T^*M)$ of the
cotangent bundle of $M$ with the pair of pants product, is
isomorphic to $H(M^{S^{1}})$ with the Chas-Sullivan product. In
their works Cohen and Jones \cite{CJ}, and Cohen, Jones and Yan
\cite{CJJ}  describe the Chas-Sullivan product in terms of a ring spectrum structure of the Thom
spectrum of a certain virtual bundle over $M^{S^{1}}.$ A most
interesting feature of this approach is that it reveals that the
essential technical point behind the Chas-Sullivan product lies in
the  construction of the so called "umkehr" map
$$F_!:  H(\mathcal{M}) \longrightarrow H(\mathcal{N})$$ for maps
$F$ between infinite dimensional manifolds under suitable
conditions, e.g. if $F$ fits into a pull-back diagram
\[\xymatrix @R=.3in  @C=.5in
{\mathcal{M} \ar[d]
\ar[r]^{F} & \mathcal{N} \ar[d]\\
M \ar[r]^{f} & N}\] where the vertical arrows are  Serre
fibrations, and $f$ is a smooth map between compact oriented
manifolds. This construction is quite general and adaptable to a
variety of context well beyond the product in string topology
\cite{ChU}. A fundamental observation by Sullivan \cite{DS} is
that in addition to the product string homology $\Ho(M^{S^{1}})$
comes with a natural co-associative co-product
$$\Ho(M^{S^{1}}) \longrightarrow \Ho(M^{S^{1}}) \otimes \Ho(M^{S^{1}}).$$
The co-product can also be explained using the Cohen-Jones
technique, indeed, in greater generality Cohen and Godin
\cite{Ch3} have constructed operations $$\mu_g:\Ho(M^{S^{1}})^{\otimes
n} \longrightarrow \Ho(M^{S^{1}})^{\otimes m}$$ on string homology
associated to each surface of genus $g$ with $n$-incoming boundary
components and $m$-outgoing boundary components. Moreover they
show that the maps $\mu_g$ give $\Ho(M^{S^{1}})$ the structure of
a topological quantum field theory in a restricted sense, i.e.
there should be a positive number of outgoing boundary components.
We remark that in a recent work \cite{tama2} Tamanoi has argued
that in most cases, e.g. if $g >0$, the operators $\mu_g$ must
vanish.\\

A natural generalization of string homology arises if one
considers the space of maps
$$M^{S^d}=\{ x \ | \ x:S^d \longrightarrow M \mbox{ smooth and
constant around the north pole} \}$$ from a $d$-dimensional sphere
into a compact oriented manifold, string topology being the case
$d$=$1.$ Thus we let $M^{S^{d}}$ be the space, with the
compact-open topology, of smooth maps from $S^{d}$ to $M$ constant
in a neighborhood of the north pole, and $C(M^{S^{d}})$ be the
graded space of chains in $M^{S^{d}}.$ In Section
\ref{sec1} we introduce the notion of transversal framed $d$-algebras,
which is based upon the notion of $d$-algebras introduced by
Kontsevich in \cite{KO1}.  After a degree shift on the complex
$C(M^{S^{d}})$ one can show the following result:\\

{\noindent}{\bf Theorem \ref{inco}.}  $\Ca(M^{S^{d}})$ is a transversal framed $d$-algebra.\\

{\noindent} Theorem \ref{inco} implies, passing to homology, a
result of Sullivan and Voronov \cite{CoVo, V1}, concerning the
algebraic
structure on the homology groups of the spaces $M^{S^{d}}$.\\

It was realized early on in string theory that alongside closed
strings it was necessary to consider open strings. A proper
understanding of open strings requires the introduction of
$D$-branes which are Dirichlet boundary conditions for the
endpoints of the open string. Perhaps the main weakness of string
theory is that actually it is not a unique theory but rather
allows for a high dimensional moduli space of models. Thus in a
sense the main open problem in the string approach towards
unification is to unify string theory itself. Various approaches
have been proposed. A promising one is the so called $M$-theory
which may be thought as a theory whose fundamental object is a
membrane moving in a given ambient manifold. This approach
stimulated the study of branes not just as boundary conditions but
as fundamental objects in their on right. In particular in
$M$-theory the dynamics of a membrane in $11$ dimensions has been
proposed as a unifying theory out of which the various models of
string theory are obtained as boundary limits. One of the main
topics of this work, developed in Section \ref{sec5}, is the study
with tools from classical algebraic topology of the space of
dynamic branes with $D$-branes as boundary conditions. Let $I$ be
the interval $[-1,1]$ and $Y$ be a compact oriented brane whose
dynamics in another manifold $M$ we like to understand. The
corresponding configuration space is
$$ M^{Y\times I }$$ the space of smooth maps from $Y\times I$ into $M$. Notice that
$Y\times I$ comes with two marked sub-manifolds, namely, its
boundary components $Y\times\{-1\}$ and $Y\times \{1\}$. For
technical reason that will become clear we restrict our attention
to a subset of possible motions for $Y$. Assume that $N_{0}\
\mbox{and}\ N_{1}$ are compact oriented embedded sub-manifolds of
$M$, then  we define the space $M^{S(Y)}(N_{0},N_{1})$ of
$Y$-branes in $M$ moving from $N_{0}$ to $N_{1}$, as follows:
\[M^{S(Y)}(N_{0},N_{1})= \left \{\gamma \ \Big |
\begin{array}{c}
\gamma \in  M^{Y\times I},\ \gamma(Y\times \{0\}) \in N_{0}, \ \gamma(Y\times \{1\}) \in
 N_{1}\\
 \gamma \mbox{ is constant around }

Y\times \{ -1 \} \mbox{ and } Y\times \{ 1 \} \\
\end{array}\right\}.
\]

By definition maps in $M^{S(Y)}(N_{0},N_{1})$ are smooth maps that
collapse the boundary components $Y\times\{-1\}$ and $Y\times
\{1\}$ to points that live in $N_0$ and $N_1$,
respectively. Once we have fixed our spaces of $Y$-branes we
construct a product for transversal pairs of chains of $Y$-branes,
i.e. we define a product
\[ \Ca(M^{S(Y)}(N_{0},N_{1})) \otimes
\Ca(M^{S(Y)}(N_{1},N_{2})) \longrightarrow
\Ca(M^{S(Y)}(N_{0},N_{2}))\] that generalizes the Sullivan product
for open strings \cite{DS} which is obtained in the case that $Y$
is a single point. This product induces well-defined product on
the corresponding homology groups
\[ \Ho(M^{S(Y)}(N_{0},N_{1})) \otimes \Ho(M^{S(Y)}(N_{1},N_{2})) \to
\Ho(M^{S(Y)}(N_{0},N_{2})),\]  which, after an appropriated degree
shift, allows us to construct a new topological invariant which
assigns to each compact oriented manifold $M$ the graded category
$\Ho(M^{S(Y)})$ whose objects are embedded oriented sub-manifolds
of $M$, and whose morphisms from $N_{0}$ to $N_{1}$ are homology
classes of $Y$-branes extended from  $N_{0}$ to $N_{1}$, i.e.
$$\Ho(M^{S(Y)})(N_{0}, N_{1})= \Ho(M^{S(Y)}(N_{0},N_{1})) .$$
Compositions are defined with the help of the product mentioned
above. We are actually going to proof a stronger result: we show
that there is a natural structure of transversal $1$-category on
the differential graded pre-category $\Ca(M^{S(Y)})$ whose objects
are embedded oriented sub-manifolds of $M$, and whose morphisms
$\Ca(M^{S(Y)}(N_{0},N_{1}))$ are chains of $Y$-branes in $M$
extended from $N_{0}$ to $N_{1}$. Moreover, after discussing some
needed notions in universal algebra such as transversal
$1$-categories, transversal traces with values in a right
$\mathcal{O}$-module where $\mathcal{O}$ is an operad, we show the
following result:\\

{\noindent}{\bf Theorem \ref{to3}.}  $\Ca(M^{S(Y)})$ is a
transversal $1$-category with  a natural $\Ca(S^{1})$-trace.\\

Section \ref{sec7} contains the main result of this work. We
introduce the notion of the homological quantum field theory HLQFT
which, in a sense, summarizes and extends the results of the
previous sections. Essentially we construct new topological
invariants for compact oriented manifold using the same basic idea
that we have been developing, but instead of considering a
correspondence of the form
$$M \longrightarrow H(M^L)$$
for fixed $L$, we consider how all this correspondences fit
together as $L$ changes. The first part of  Section
\ref{sec7} may be regarded as a second introduction to this work and presents
a general panorama of the categorical approach to the definition
of quantum field theories. Let us here just highlight the main
ingredients involved in our notion.   The main object is the
theory of cobordisms introduced by Thom in \cite{RT, RT1}. From a
physical point of view we may think of the theory of cobordisms as
the theory of space and their interactions trough space-time. A
subtle but fundamental issue is that both space and space-time may
be disconnected. Another delicate issue that the empty space has
to be included as a valid one. Using Thom's cobordisms Atiyah
\cite{MA} wrote down the axioms for topological quantum field
theory TQFT, a type of quantum field theory that had been
introduced earlier by Witten in \cite{W1, WZ}. TQFT are  of great
interest for mathematicians since the  vacuum to vacuum
correlation functions of such theories are by construction
topological invariants for compact oriented manifolds. The
Atiyah's axioms for TQFT essentially (omitting unitarity) identify
the category of topological quantum field theories with the
category of monoidal functors from the category of cobordisms into
the category of finite dimensional vector spaces.  A further
development in the field was the introduction by Turaev in
\cite{TU1, TU2} of homotopy quantum field theories, following a
pattern similar to the one explained above for TQFT, but replacing
the category of cobordisms by a certain category of cobordisms
provided with homotopy classes of maps into a given topological
space. A homotopy quantum field is a monoidal functor from that
generalized category of cobordisms into vector spaces. In order to
define homological quantum field theories we first introduce the
notion of cobordisms provided with homology classes of maps into a
fixed compact oriented smooth manifold. Next, we defined a HLQFT
as a monoidal functor from that category of extended cobordisms
into the category of vector spaces. In contrast with Turaev's
definition, we demand that the maps from cobordisms to the fixed
manifold be constant on a neighborhood of each boundary component.
This is a major technical restriction which is necessary in order
to define composition of morphisms using transversal intersection
on finite dimensional manifolds. Without imposing this restriction
one is forced to deal with transversal intersections on infinite
dimensional manifolds, a rather technical subject that we prefer
to avoid in this paper. In Sections \ref{sss} and \ref{ttt} we
give examples and discuss the possible applications  of
homological quantum field theories in dimensions $1$ and $2$,
respectively.

\subsection*{Acknowledgment}
Our thanks to Raymundo Popper for his support and interest in the
ideas of this work. The first author thanks Alexander Cardona and
Sylvie Paycha for inviting him to deliver a talk  at the school
"Geometric and Topological Methods for Quantum Field Theory,"
Colombia 2007. The second author thanks Nicolas Andruskiewisth for
inviting him to participate in the "XIV Coloquio Latinoamericano
de Algebra," Argentina 2001, where he attended lectures on "String
Topology" by Dennis Sullivan. He thanks MSRI and the organizers of
the conference "String Topology,"  Mexico 2006, for inviting him
to participate. He also thanks Takashi Kimura and Bernardo Uribe
for inviting him to deliver talks at Boston University and at the
Universidad de los Andes, respectively.

\section{Transversal algebras and categories}\label{sec1}

In this section we first introduce the basic background needed to
state and prove the main results of this work. There are two
fundamental ingredients that we shall need:
\begin{itemize}
\item We must be able to work with algebras, and categories, with products
defined only for transversal sequences.

\item We need to introduce an appropriated chain model for the homology of smooth manifolds
such that the transversal intersection of chains becomes a
transversal algebra.
\end{itemize}
Once we are done with these preliminary constructions, we apply
them to study the algebraic structure on the space $C(M^{S^{d}})$
of chains of maps form the $d$-sphere into a compact oriented
manifold $M$.\\

In this work all vector spaces are defined over the complex
numbers. We denote by $\de$ the symmetric monoidal category of
differential $\mathbb{Z}$-graded vector spaces. Objects in $\de$
are pairs $(V,d),$ where
$$V=\bigoplus_{i\in \mathbb{Z}}V_{i}$$ is $\mathbb{Z}$-graded
vector space and $d\colon V \longrightarrow V$ is such that
$d_{i}\colon V_{i}\longrightarrow V_{i-1}$ and $d^{2}=0.$  A
morphism $\f:(V_1,d_1) \longrightarrow (V_2,d_2)$ is a degree
preserving linear map $f:V_1 \longrightarrow V_2$ such that
$d_2f=fd_1. $ For each $n \in \mathbb{Z}$, right tensor
multiplication with the complex  $\mathbb{C}[n]$ such that
$\mathbb{C}[n]^{-n}=\mathbb{C}$ and $\mathbb{C}[n]^{i}=0$ for
$i\neq -n,$ gives a  functor
$$[n]:dg\mbox{-}vect \longrightarrow dg\mbox{-}vect $$
which sends $V$ into $V[n]$ where $V[n]^{i}=V^{i+n}.$ We say that
$V[n]$ is equal to $V$ with degrees shifted down by $n.$ To
simplify notation at various stages in this work in which a degree
shift is fixed within a given context, we shall write $\mathrm{V}$
for the vector space with shifted degrees. For example if a shift
of degree by $n$ is involved then $\mathrm{V}=V[n].$ A
differential graded precategory or dg-precategory $\C$ consists of
following data:
\begin{itemize}
\item A collection of objects $\ob(\C).$

\item For $x,y \in \ob(\C)$ a differential graded vector space $\C(x,y)$
 called the space of morphisms from $x$ to $y.$
\end{itemize}
A prefunctor $\f\colon \C \longrightarrow \d$ consists of a map
$\f\colon \ob(\C) \longrightarrow \ob(\d)$ and for each pair of
objects $x,y \in \ob(\C)$ a morphism of differential graded vector
spaces
$$\f_{x,y}\colon\C(x,y)\longrightarrow \d(\f(x),\f(y)).$$
We define graded precategories or g-precategories as
dg-precategories with vanishing differentials on the spaces of
morphisms.  The homology $H(\C)$ of a dg-precategory $\C$ is the
g-precategory given by:
\begin{itemize}
\item $\ob(H(\C))=\ob(\C)$.

\item $H(\C)(x,y)=H(\C(x,y))$  for objects $x,y$
of $H(\C).$
\end{itemize}

Notice that if $\C$ is actually a category, i.e. in addition to
the structure of pre-category it has compositions and identities,
then $H(\C)$ is also a category with the induced composition maps.
We do not want to restrict ourselves to consider only the case
where $\C$ is a category for two reasons. On the one hand, we
shall consider more general structures than simple categories, for
example, structures where there are not just one but a whole set
of different ways to compose morphisms. On the other hand, we are
interested in the case where the compositions of morphisms in $\C$
are not a quite defined for all morphisms,  but only for some sort
of distinguished sequences of morphisms called transversal
sequences. So what we need is to specify the conditions for the
domain of definition of these partially defined compositions. We
assume that for each sequence of objects $x_0,...,x_n$ we have a
subspace $\C(x_{0},\cdots ,x_{n})$ of
$\bigotimes_{i=1}^{n}\C(x_{i-1},x_{i})$ consisting of transversal
$n$-tuples of morphisms of $\C$, i.e. generic sequences for which
compositions are well defined.  We shall demand that any $0$ or
$1$-tuple of morphisms is automatically transversal, and that for
$n\geq 2$ any closed $n$-tuple of morphisms in
$\bigotimes_{i=1}^{n}\C(x_{i-1},x_{i})$ is homologous to a closed
transversal $n$-tuple. Finally we demand that any subsequence of a
transversal sequence be transversal. We formalize these ideas in
our next the definition which is modelled on the corresponding
notion for algebras given by Kriz and May  \cite{KM}.

\begin{defi}\label{d1}
{\em A domain $\C^{*}$ in a dg-precategory $\C$ consists of the
following data:}

\begin{enumerate}

\item {\em  A differential graded vector space $\C(x_{0},\cdots
 ,x_{n})$ for $x_{0},\cdots ,x_{n}$ objects of $\C.$}

\item  {\em $\C(\emptyset)=k.$}

\item {\em Inclusion maps $i_{n}\colon
\C(x_{0},\cdots ,x_{n}) \longrightarrow
\bigotimes_{i=1}^{n}\C(x_{i-1},x_{i}).$}
\end{enumerate}

{\em This data should satisfy the following properties:}
\begin{enumerate}

\item {\em $i_{\emptyset}: k \longrightarrow k$ is the identity map.}

\item {\em $i_{1}\colon \C(x_{0},x_{1}) \longrightarrow
\C(x_{0},x_{1})$ is the identity map.}

\item {\em $i_{n}\colon \C(x_{0},\cdots
,x_{n}) \longrightarrow
\bigotimes_{i=1}^{n}\C(x_{i-1},x_{i})$ is a quasi-isomorphism,
i.e. the induced map
$$H({i}_{n})\colon H(\C( x_{0},\cdots ,x_{n})) \longrightarrow
\bigotimes_{i=1}^{n} H(\C(x_{i-1},x_{i}))$$ is an
isomorphism.}

\item {\em  For a partition $n=n_{1}+ \cdots + n_{k}$ of $n$ in $k$ parts, we set
$m_{0}=0$ and $m_{i}=n_{1}+\cdots +n_{i},$ for $1\leq i \leq k.$
The inclusion map $i_{n}$  should factor through
$$\bigotimes_{i=1}^{k}\C(x_{m_{i-1}},\cdots
 ,x_{m_{i}})$$
 as indicated in the following commutative diagram}
\[ \xymatrix@R=.4in  @C=.3in{ \des
\bigotimes_{i=1}^{k}\C(x_{m_{i-1}},\cdots
 ,x_{m_{i}}) \ar[d]_{\des\otimes i_{n_{i}}}&
\C(x_{0},\cdots ,x_{n}) \ar[d]^{i_{n}} \ar@{_{(}->}[l] \\
\des\bigotimes_{i=1}^{k}
\bigotimes_{s=m_{i-1}+1}^{m_{i}}\C(x_{s-1},x_{s}) \ar[r]_{\ \
\ \ =}  &
\des\bigotimes_{i=1}^{n}\C(x_{i-1},x_{i})} \]
\end{enumerate}
\end{defi}

In order to formally introduce the possibility of multiple types
of compositions, we need to recall the notion of operads defined
in a symmetric monoidal category with product $\otimes$ and unit
object $1$; typical examples of latter kind of  categories, and
the only ones that will be consider in this work, are the
categories of sets, topological spaces, vector spaces, graded
vector spaces and differential graded vector spaces. A
non-symmetric operad $\O$ consists of a sequence $\O_n,$ for $n
\geq 0$, of objects in the corresponding category, an unit map
$\eta
\colon 1 \to\O_1,$  and maps
\[ \gamma_{k} \colon \O_k \otimes \O_{n_{1}} \otimes \cdots
\otimes \O_{n_{k}} \longrightarrow \O_{n_1+...+n_k} \] for $k\geq 1$ and $n_{s}\geq 0.$
The maps $\gamma_{k}$  are required to be associative and unital
in the appropriated sense. The reader will find a lot information
about operads in
\cite{May3}, see also \cite{get1} for a recent fresh approach. If in addition a right action of the symmetric
group $S_{n}$ on $\O_n$ is given and the maps $\gamma_{n}$ are
equivariant, then we say that $\O$ is an operad. To any object $x$
in a symmetric monoidal category, there is attached an operad,
called the endomorphisms operad, with is $n$ component given by
$$End_x(n)=Hom(x^n,x)$$ For a given operad $\O$ in the same
category, one says that $x$ is a $\O$ algebra, if there is a
morphisms of operads $\theta:\O \longrightarrow End_x,$ i.e. a
sequence of maps $\theta_k:\O \otimes x^k \longrightarrow x$
satisfying certain natural axioms. It is easy to check that there
are operads whose algebras are exactly associative algebras,
commutative algebra, Lie algebras, Poisson algebra, BV algebras,
$A_{\infty}$-algebras, $A_{\infty}^N$-algebras
\cite{angel2}, etc. One can in a similar fashion define for each
operad the category $\O$-categories. We shall not make explicit
that definition since we are presently going to consider the more
general notion of partial $\O$ categories.

\begin{defi}\label{d2}
 {\em Let $\O$ be a non-symmetric dg-operad and $\C$ be a dg-precategory.
We say that $\C$ is a transversal $\O$-category if the following
data is given:}
\begin{enumerate}

\item  {\em A domain $\C^{*}$ in $\C.$}

\item   {\em Maps $\theta_{n}\colon \O_n \otimes \C(x_{0},\cdots
,x_{n}) \longrightarrow \C(x_{0},x_{n})$ for  $x_{0}, \cdots
,x_{n}
\in
\ob( \C).$}
\end{enumerate}

  {\em This data should satisfy the following axioms:}

\begin{enumerate}
\item  {\em $\theta_{1}(1\otimes
\C(x_{0},x_{1}))=\C(x_{0},x_{1})$ where $1$ denotes the identity
in $\O(1).$}

\item  {\em For $n=n_{1}+ \cdots + n_{k},$ set
$m_{0}=0$ and $m_{i}=n_{1}+\cdots +n_{i}.$ The maps
$$\des\bigotimes_{i=1}^{k}\O_{n_{i}}\otimes
\des\bigotimes_{i=1}^{k} \C(x_{m_{i-1}},\cdots
,x_{m_{i}})\longrightarrow
\des\bigotimes_{i=1}^{k} \C(x_{m_{i-1}} ,x_{m_{i}})$$

obtained by the composition of the inclusion  $$\C(x_{0},\cdots
,x_{n}) \subseteq \bigotimes_{i=1}^{k}
\C(x_{m_{i-1}},\cdots ,x_{m_{i}}),$$ shuffling, and the application
of  $\theta ^{\otimes k},$ factors through
$\C(x_{m_{0}},\cdots,x_{m_{k}})$ as indicated in the following
diagram:}
\[\xymatrix@R=.3in @C=1in{\des
\bigotimes_{i=1}^{k}\O_{n_{i}}\otimes \C(x_{0},\cdots ,x_{n})
\ar[r] \ar@ {.>}[d] & \des\bigotimes_{i=1}^{k}
\C(x_{m_{i-1}} ,x_{m_{i}})\\
\C(x_{m_{0}},\cdots ,x_{m_{k}})\ar @
{^{(}->}[ru]  &  }\]

\item {\em The following  diagram commutes}
\[\!\!\!\!\!\!\!\!\!\!\!\!\!\!\!\!\!\! \xymatrix @R=.3in
@C=.2in {\O_k \otimes \des
\bigotimes^{k}_{s=1}\O_{n_{s}} \otimes \C(x_{0},\cdots ,x_{n})
\ar[r]^{\ \ \ \ \ \ \gamma_{k} \otimes \id} \ar @ {^{(}->}[d]\ar@
{.>}[rd]  & \O_{n}\otimes \C(x_{0},\cdots ,x_{n})
\ar[rd]^{\theta_{n}} & \\
\O_k\otimes \des
\bigotimes^{k}_{s=1}\O_{n_{s}}
\otimes \des\bigotimes_{s=1}^{k}\C(x_{m_{s-1}},\cdots
,x_{m_{s}}) \ar[d]_{\mbox{shuffle}}
 &   \O_k\otimes \C(x_{m_{0}},\cdots,x_{m_{k}})
 \ar @ {^{(}->}[d]^{\rm{1}  \otimes \it{i}_{k} }
\ar[r]_{\ \ \ \ \ \ \ \ \ \  \theta_{k}} &
\C(x_{0},x_{n})\\
\O_k \otimes \des \bigotimes^{k}_{s=1}
\O_{n_{s}} \otimes \des\bigotimes_{s=1}^{k} \C(x_{m_{s-1}},\cdots
,x_{m_{s}}) \ar[r]  & \O_k\otimes \des
\bigotimes_{i=1}^{k} \C(x_{m_{i-1}} ,x_{m_{i}}) & }\]
\end{enumerate}
\end{defi}

Notice that in the definition above $\O$ is an operad in the
standard sense, i.e. compositions are always well defined at the
operadic level. What is transversally defined is the action of the
operad $\O$ on the precategory $\C$, i.e. the various composition
of morphisms in $\C$. In our next result we use the known fact
that if $\O$ is a dg-operad, then the sequence $H(\O)$ given by
$H(\O)_n=H(\O_n)$ with the structural maps induced from those of
$\O$ is a g-operad.

\begin{thm}\label{Teo4}
 {\em If $\C$ is a transversal $\O$-category then $H(\C)$ is an
$H(\O)$-category.}
\end{thm}

Indeed if $\C$ is a transversal $\O$-category then the morphism
$$i_{n}\colon
\C(x_{0},\cdots ,x_{n})\longrightarrow \des\bigotimes_{i=1}^{n}
\C(x_{i-1} ,x_{i})$$ is a quasi-isomorphism. Hence the vertical
arrow in the  diagram
\[\xymatrix@R=.2in  @C=.4in {H(\O_n)\otimes
H(\C(x_{0},\cdots,x_{n})) \ar[r] \ar[d] & H(\C(x_{0},x_{n}))\\
H(\O_n)\otimes
\des\bigotimes_{i=1}^{n}H( \C(x_{i-1} ,x_{i})) \ar[ru] & }\]
is also an isomorphism. The diagonal map above gives  $H(\C)$ the
structure of  a $H(\O)$-algebra.\\

The concept of transversal $\O$-algebra where $\O$ is a
non-symmetric operad is easily deduced from that of a transversal
$\O$-category. Indeed we say that  a differential graded vector
space $A$ is a transversal $\O$-algebra if the precategory
$\C_{A}$ with a unique object  $p$ such that $\C_{A}(p,p)=A$ is an
$\O$-category. If $\O$ is an operad we demand in addition that the
maps $\O_n\otimes A_{n}\longrightarrow A$ be $S_{n}$-equivariant.
The reader may consult \cite{KM} where various interesting
examples of transversal $\O$-algebras are studied. Another
interesting example was introduced by Karoubi in
\cite{MAX, MAX1} where he associated to each simplicial set $X$ its transversal $\mathbb{Z}$-algebra of
quasi-commutative cochains, which determines the homotopy type of
$X$.

\begin{cor}\label{cor5}
 {\em If $A$ is a transversal $\O$-algebra, then $H(A)$ is an
$H(\O)$-algebra.}
\end{cor}

This closes our comments on the structure of transversal
categories. Let us next consider the chain model that we are going
to be using along this work. In a nutshell, for a given manifold
$M$, our space of chains in $M$ is generated by smooth maps from
manifolds with corners into $M$. We recall that a $n$-dimensional
manifold with corners $M$ is naturally a stratified manifold
$$M=
\des \bigsqcup_{0\leq l \leq n} \partial _{l} M$$ where
the smooth strata are the connected components of $\partial _{l}
M$ where
\[\partial_{l} M=\{m\in M\mid \ \mbox{there exists local
coordinates mapping } m
 \mbox{ to} \ \partial_{l}H^{n}_{k} \}\]
where $$H^{n}_{k}=[0,\infty)^{k}\times
\mathbb{R}^{n-k}=\{(x_{1},\cdots ,x_{n})\in \mathbb{R}^{n} \ \ | \
\ x_{1} \geq 0, \cdots ,x_{k} \geq 0 \}$$ and $$\partial_{l}H^{n}_{k}=\{x\in H^{n}_{k}\mid x_{i}=0 \
\mbox{for exactly}\ l \ \mbox{of the first}\ k \ \mbox{indices}\}.$$
Given an oriented manifold $M$ we define the graded vector space
\[C(M)=\bigoplus_{i\in \N} C_{i}(M),\]
where $C_{\it i}(M)$   denotes the complex vector space
constructed as follows:
\begin{itemize}
\item{Let $\overline{C}_{i}(M)$ be the vector space freely generated by equivalence classes of pairs
$(K,c)$ where $K$ is a compact oriented manifold with corners and
$c:K
\longrightarrow M$ is a smooth map.  A pair $(K,c)$ is equivalent to
another $(L,d)$ if and only if there exists a orientation
preserving diffeormorphism $f:K \longrightarrow L$ such that $d
\circ f =c.$  Abusing notation the equivalence class of $(K,c)$ is also denoted by
$(K,c)$. The collection of equivalence classes of such pairs is a
set since any manifold with corners is diffeomorphic to a manifold
with corners embedded in some $\mathbb{R}^n$, and thus one can
assume that the domain $K$ of all chains are embedded in
$\mathbb{R}^n$ for some $n \in \mathbb{N}$.}

\item{$C_{i}(M)$ is the quotient of $\overline{C}_{i}(M)$ by the
following relationships:\\

1)$(K^{op},c) = -(K,c)$ where $K^{op}$ is the manifold $K$
provided with the opposite orientation.\\

2) $(K_1 \sqcup K_2, c_1 \sqcup c_2)= (K_1,c_1) +(K_2,c_2)$.}
\end{itemize}

\begin{figure}[h]
\centering
\includegraphics[width=0.5\textwidth]{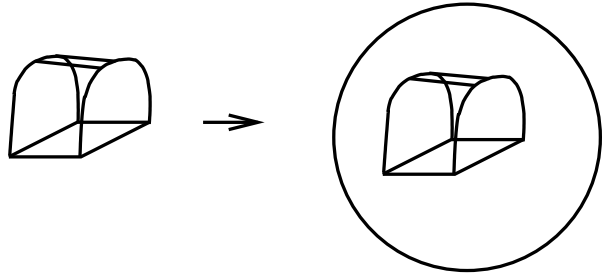} \caption{ \ Example of a chain with domain a manifold with
corners.}\label{Bor77}
\end{figure}

Figure \ref{Bor77} shows an example of a chain with domain a
manifold with corners.

We define a  differential $\partial \colon C_{i}(M)
\longrightarrow C_{{i}-\rm {1}}(M)$ on $C(M)$ as follows:
$$\partial(K_,c)=\sum_{L \in \pi_{0}(\partial_{1}
K)}(\overline{L},c|_{ \overline{L}}),$$ where the sum ranges over
the connected components of the first boundary strata
$\partial_{1} K$ of $K$ provided with the induced boundary
orientation. We denote by $c|_{\overline{L}}$ the restriction of
$c$ to the closure of $L.$ Complexes $C(M)$ enjoy the following
crucial property that shows that we can  compute singular homology
using the manifold with corners chain model.

\begin{thm}
 {\em $(C(M),\partial)$ is a differential $\mathbb{Z}$-graded
vector space. Moreover}
\[H(\Ca(M),\partial)=H(M)=\mbox{\em singular homology of }
M.\]
\end{thm}

In fact the identity
$$\partial^{2}(K,c)= \sum_{L \in \pi_0(\partial_{2} K) }[(\overline{L},
c|_{\overline{L}}) + (\overline{L}^{op}, c|_{\overline{L}})] $$
implies that $\partial^{2}=0.$ There is an obvious inclusion
$i\colon C^{s}(M)\longrightarrow C(M)$  of the complex of singular
chains into the complex of chains with manifolds with corners as
domain of definitions. The map $i$ is a quasi-isomorphism since
any manifold with corners can be triangulated \footnote{We thank
J. Brasselet, M. Goresky, R. Melrose, and A. Dimca for helpful
comments on the triangulation of manifold with corners. See
references
\cite{MG, JEF, ve} for more information on triangulability.} and
thus any chain in $C(M)$ is homologous to a chain in $C^s(M).$\\

\begin{figure}[h]
\centering
\includegraphics[width=0.4\textwidth]{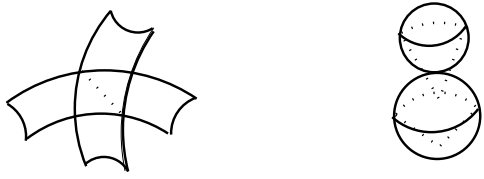} \caption{ \ Left: transversal submanifolds. Right: non-transversal
submanifolds.}\label{Bor56}
\end{figure}

The definition as well as many results for transversal smooth maps
can be generalized along the lines of \cite{gui, guill} so that
they apply to maps from manifolds with corners into smooth
manifolds. Recall that two submanifolds $K$ and $L$ of a smooth
manifold $M$ are transversal if for each $x \in K
\cap L$ one has that:
$$T_xK + T_xL = T_xM.$$
The remarkable fact is that if $K$ and $L$ are transversal, then
$K\cap L$ is also a submanifold of $M$. Figure \ref{Bor56} shows a
transversal pair, and a non-transversal pair of submanifolds of
$\mathbb{R}^3.$

 Likewise if $f\colon K
\longrightarrow M$  is a smooth map from a
manifold with corners and $L$ a submanifold of $M$, then we say
that $f$ is transversal to $L$ if for $x \in \partial_s K$ such
that $f(x) \in L$ we have that
$$df(T_x\partial_s K) + T_{f(x)}L=T_{f(x)}M .$$
One can check that in this situation then the pre-image
$f^{-1}(L)$ is a submanifold with corners of $K$, and that the
co-dimension of $f^{-1}(L)$ is equal to the co-dimension of $L$.
The notation $f \pitchfork L$ means that the map $f$ is
transversal to the submanifold $L$. Next, assume that we have maps
$f_1:K_1
\longrightarrow M,..., f_n:K_n
\longrightarrow M$ from manifolds with corners into $M$. We say
that the maps $f_1,...,f_n$ are transversal if the map
$$(f_1,...,f_n):K_1 \times ... \times K_n \longrightarrow M ,$$
is transversal to $\Delta_{n}$, the $n$-diagonal submanifold of
$M^{n}$ given by
$$\Delta_{n}=\{ (m,\ldots ,m) \ | \ m \in M \}.$$
In this case we have that
$$(f_1,...,f_n)^{-1}(L)=K_1 \times_{M} ... \times_{M}K_n$$
is a submanifold with corners of $K_1 \times ... \times K_n$ of
co-dimension $(n-1)\dim M.$ For example if $K$ and $L$ are
manifolds with corners, $f\colon K
\longrightarrow M$ and $g\colon L
\longrightarrow M$ are smooth maps. Then $f$ and  $g$ are
transversal maps if for  $0\leq s\leq \dim K, \ 0\leq t\leq \dim
L$ the restrictions of $f$ and $g$ to $\partial_{s}K$ and
$\partial_{t}L$, respectively, are transversal maps, i.e. given $x
\in\partial_{k}K$ and $y\in
\partial_{s}L$ such that  $f(x)=g(y)=m$,  we must have that
\[df(T_{x}\partial_{k}K)
+dg(T_{y}\partial_{s}L)=T_{m}M.\] In short, two maps are
transversal if their respective restrictions
to the smooth strata are transversal. \\

One of the main advantages of the category of manifolds with
corners is that unlike the category of manifolds with boundaries
it is closed under Cartesian products, and even more remarkably it
is generically closed under fibred products. Indeed with the
notion of transversality given above one can show the following
result \cite{cd}. Let $ K_{x}$, $K_{y}$ and $K_z$ be oriented
manifolds with corners and $M$ be an oriented smooth manifold.

\begin{thm}\label{l1}
 {\em Let $x\colon K_{x}\longrightarrow M,$ $y\colon K_{y}\longrightarrow M$ and $z\colon K_{z}\longrightarrow M$ be
transversal smooth maps, then
\begin{itemize}
\item $K_{x}\times _{M}K_{y}= \{ (a,b)\in
K_{x}\times K_{y}\mid x(a)=y(b) \}$ is in a natural way an
oriented manifold with corners embedded in $K_{x}\times K_{y}$.

\item $(K_{x}\times_{M}K_{y})\times_M K_z = K_{x}\times_{M}K_{y}\times_M
K_z= K_{x}\times_{M}(K_{y}\times_M K_z).$

\item $K_x \times_M K_y= (-1)^{(\dim K_x + \dim M)(\dim K_x + \dim M)}K_y \times_M
K_x.$

\end{itemize}}
\end{thm}

We are ready to study the algebraic structure on the space of
chains of maps from spheres into a given compact oriented manifold
$M$. We  let $$D^{d}=\{ x\in
\mathbb{R}^{d}
\mid x_1^2 + ... + x_d^2=1
\}$$ be the unit disc in $\mathbb{R}^{d}.$ By definition a little disc in $D^{d}$
is an affine transformation $$T_{a,r}\colon D^{d}\longrightarrow
D^{d}$$ given by $T_{a,r}(x)=rx+a,$  where $0 < r < 1$ and $a \in
D^{d}$ are such that $\im (T_{a,r}) \subseteq D^{d}.$  For $n\geq
0,$ consider the spaces
\[D^{d}_n= \Big\{(T_{a_{1},r_{1}}, \ldots ,T_{a_{n},r_{n}}) \ \Big |
\begin{array}{c}
a_{i}, a_{j}\in D^{d},\ 0<  r_{i} < 1 \ \mbox{such that if }\
i\neq j\
 \\
\mbox{then} \ \ \overline{\im(T_{a_{i},r_{i}})}\cap \overline{\im(T_{a_{j},r_{j}})}=
\emptyset  \\
\end{array} \Big \} .\]

\begin{figure}[h]
\centering
\includegraphics[width=0.7\textwidth]{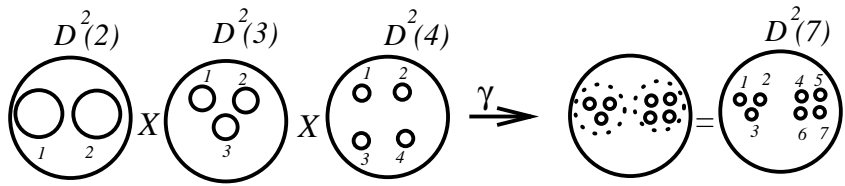} \caption{ \ Composition in the operad of
little discs in the plane.}\label{Bor14}
\end{figure}

Notice that the disc with center $a$ and radius $r$ is obtained as
the image of the transformation $T_{a,r}$ applied the standard
disc $D^{d}$. The sequence of topological spaces $D^{d}_n$ carries
a  natural structure of operad, called the little $d$-discs operad
and denoted by $D^{d}$. The little disc operad was introduced by
Boardman and Vogt, in its cubic version, in
\cite{BV} and May  in \cite{May1}. Figure
\ref{Bor14} illustrates how compositions are defined in the operad of
little discs.

The framed little $d$-discs operad $fD^{d}$ is obtained by placing
an element of $SO(d)$, the group of orientation and metric
preserving linear transformations of Euclidean $d$-space, on each
little disc. Explicitly we have that
$$fD^{d}_n=D^d_n \times SO(d)^n. $$ The composition
\[ \gamma_{k} \colon fD^{d}_k \times fD^{d}_{n_{1}} \times \cdots
\times fD^{d}_{n_{k}} \longrightarrow fD^{d}_{n_1+...+n_k} \]
is given by $$\gamma_{k}[(c,g),(b_1,h^1),...,(b_k,h^k)]=
(\gamma_{k}(c,g_1b_1,...,g_kb_k),g_1h^1,...,g_kh^k),$$ where
$h^i=(h_1^i,...,h_{n_i}^i)$ and
$g_ih^i=(g_ih_1^i,...,g_ih_{n_i}^i).$\\

Notice that the defining action of $SO(d)$ on $\mathbb{R}^d$
induces an action $SO(d)\times D^d \longrightarrow D^d$ which in
turns induces actions $ SO(d)\times D^d_n \longrightarrow D^d_n$.
It is not hard to see that an algebra $X$ over $fD^{d}$ is the
same as an algebra $X$ over $D^d$ provided with a  $SO(d)$ action
$SO(d)\times X \longrightarrow X$ such that
$$\theta_k(a,g_1x_1,...,g_kx_k)=g\theta_k(a,x_1,...,x_k)  .$$
Moreover in this case the $fD^{d}$-structure on $X$ is related to
the $fD^{d}$-structure on $X$ as follows:
$$\theta_{fD^d}((a,g_1,...,g_k),x_1,...,x_k)=\theta_{D^d}(a,g_1x_1,...,g_kx_k) .$$

\begin{defi}
 {\em A  transversal framed $d$-algebra is a transversal algebra over the operad $C(fD^{d})$ of
 chains in the framed little $d$-discs  operad.}
\end{defi}

\begin{defi}\label{def13}
 {\em Let $M^{S^{d}}$ be the set of  smooth maps $\alpha
\colon D^{d}\longrightarrow M$  constant in an open
neighborhood of $\partial( D^{d}).$ We topologize $M^{S^{d}}$ with
the compact-open topology.}
\end{defi}

Next we have to deal with a rather subtle and fundamental issue.
We like to study the homology and more generally the chains on the
space $M^{S^{d}}.$ Above we introduced a chain model for smooth
manifolds, where a chain is a smooth map from a manifold with
corners into the manifold in question. Of course $M^{S^{d}}$ is
not a manifold in the usual sense since it is an infinite
dimensional space. However with can avoid running into troubles by
adopting the following convenient definition for the space of
chains in $M^{S^{d}};$ it is straightforward to check that with
this definition we obtain a chain model that indeed computes the
homology of $M^{S^{d}}.$ Thus we shall consider the vector space
$$C(M^{S^{d}})=\bigoplus_{i=0}^{\infty}
C_{i}(M^{S^{d}})$$  generated by equivalence classes  of  maps
$x\colon K_{x} \longrightarrow M^{S^{d}}$ such that the associated
map $$\widehat{x}\colon K_{x}\times D^{d}
\longrightarrow M$$ given by $\widehat{x}(c,p)=x(c)(p)$  is a smooth map. \\

Let $e\colon M^{S^{d}}\longrightarrow M$ be the map given by
$e(\alpha)=\alpha(\partial (D^{d})).$ We shall also denote by $e$
the induced map $e\colon \Ca(M^{S^{d}}) \longrightarrow
\Ca(M)$ given by $$e(\sum a_{x}x)=\sum a_{x}e(x).$$
Given chains $x_{i}\colon K_{x_{i}}\longrightarrow M^{S^{d}}$ for
$1\leq i \leq n$  consider the map
\[\xymatrix @R=.08in  { e(x_{1},
 \ldots , x_{n})\colon &\!\!\!\!\!\!\!\!\!\!\!\!\!\!\!\!\!\!\!\!
\des\prod^{n}_{i=1} K_{x_{i}}
\ar[r] &  M^{\times n} \\
&  \!\!\!\!\!\!\!\!\!\!\!\!\!\!\!\! (c_{1}, \ldots , c_{n})\ar
@{|->}[r] &\ \ (e(x_{1}(c_{1})) , \ldots ,e(x_{n}(c_{n})) ) } \]
The map $e$ is  smooth and thus according to Theorem
\ref{l1} if $e(x_{1}, \ldots , x_{n})\pitchfork\Delta_{n}$ then
$$e^{-1}(\Delta_{n})=\{ (c_{1}, \ldots ,c_{n}) \in
\des\prod^{n}_{i=1} K_{x_{i}} \mid e(x_{1}(c_{1}))= \cdots
=e(x_{n}(c_{n}))\}$$ is a manifold with corners.\\

Consider the sequence $\Ca(M^{S^{d}})^{*}$ in dg-vect where for $n
\geq 0$ we let
$$\Ca(M^{S^{d}})^{\it n}\subseteq
\Ca(M^{S^{d}})^{\otimes n}$$ be the subspace generated
by tuples $x_{1}\otimes  \cdots \otimes x_{n}$, with $x_{i}\in
\Ca(M^{S^{d}})$ for $1\leq i\leq n$, such that
$$e(x_{1}, \cdots ,x_{n})\pitchfork \Delta_{n}.$$

\begin{thm}\label{lem 16} {\em $\Ca(M^{S^{d}})^{*}$ is
a domain in $\Ca(M^{S^{d}}).$}
\end{thm}

We check that the axioms of Definition \ref{d1} hold. Since
$\Ca(M^{S^{d}})^{0}=
\Ca(M^{S^{d}})^{\otimes 0}=k$  axiom 1 holds by convention.
Clearly $$\Ca(M^{S^{d}})^{ 1}= \Ca(M^{S^{d}})^{\otimes 1}=
\Ca(M^{S^{d}})$$ since for any chain $x\colon K_{x}\longrightarrow M^{S^{d}}$
in $\Ca(M^{S^{d}})$ one checks that $e(x)\pitchfork \Delta_{1}$ as
follows $$ \im(de(x))+\im(d\Delta_{1})\supseteq
\im(d\Delta_{1})=TM.$$ Thus axiom 2 also holds.
By Sard's lemma any chain $$(x_{1},\cdots ,x_{n})\colon
\des\prod^{n}_{i=1} K_{x_{i}}\longrightarrow (M^{S^{d}})^{n}$$ is
homologous to a chain $$(y_{1},\cdots ,y_{n})\colon
\des\prod^{n}_{i=1} K_{y_{i}}\longrightarrow (M^{S^{d}})^{n}$$ such
that $e(y_{1},\cdots ,y_{n})\pitchfork \Delta_{n}.$  Thus the
inclusion maps $i_{n}\colon
\Ca(M^{S^{d}})^{ n}\longrightarrow
\Ca(M^{S^{d}})^{\otimes n}$ are quasi-isomor\-phisms and axiom
3 holds. Axiom 4 is an obvious consequence of the definition of $\Ca(M^{S^{d}})^{\it n}$ given above.\\

Our next result provides a natural algebraic structure on the
space $$\Ca(M^{S^{d}})=C(M^{S^{d}})[\dim M]$$ of chains of maps
from the $d$-sphere into $M$ with degrees shifted down by $\dim
M.$ Notice that the action of $SO(d)$ on $D^d$ induces an action
$SO(d)\times M^{S^d} \longrightarrow M^{S^d}$.

\begin{thm}\label{inco}
{\em The dg-vect $\Ca(M^{S^{d}})$  has a natural structure of
transversal framed $d$-algebra.}
\end{thm}

In order to prove this result we must define for each $n \geq 0$ a
map
$$\theta_{n}\colon C(fD^{d}_n)
\otimes
\Ca(M^{S^{d}})^{n} \longrightarrow \Ca(M^{S^{d}}).$$ This is done as follows.  Given $x \in
C(D^{d}_n)$ and $ x_{i} \in \Ca(M^{S^{d}})$ for $1\leq i \leq n,$
then the  domain of $\theta_{n}(x;x_{1},\ldots ,x_{n})$ is the
manifold with corners given by $$K_{\theta_{n}(x;x_{1},\ldots
,x_{n})}=K_{x}
\times e(x_1,...,x_n)^{-1}(\Delta_{n} ).$$  For $c\in K_{x}$ we let $x(c)$ be
given by
$$x(c)=(T_{p_{1}(c),r_{1}(c)},\cdots , T_{p_{n}(c),r_{n}(c)}, g_1(c),...,g_n(c)).$$ The
map
\[\xymatrix @R=.08in {\theta_{n}(x;x_{1},\ldots ,x_{n})\colon
K_{\theta_{n}(x;x_{1},\ldots ,x_{n})} \ar[rr]& &\Ca( M^{S^{d}})
}\] is such that for $(c;c_{1}\cdots ,c_{n})\in
K_{\theta_{n}(x;x_{1}\cdots ,x_{n})}$ and $y\in D^{d}$ we have
that
\[ \theta_{n}(x;x_{1},\ldots ,x_{n})(c;c_{1},\ldots ,c_{n})(y) =
 \left\{ \begin{array}{ll}
e(x_{1}(c_{1}))  \ \ \mbox{if}\ y
\notin \bigcup \overline{\im(T_{p_{i}(c),r_{i}(c)})}\\
  \ &  \\
g_i(c_i)x_{i}(c_{i})\Big(\frac{y-p_{i}(c)}{r_{i}(c)}\Big) \ \
\mbox{if}\ y
\in \overline{\im(T_{p_{i}(c),r_{i}(c)})}
\end{array} \right. \]
We check axioms 1, 2 and 3 of Definition \ref{d2}. We need to
check that $\theta_{1}(1;x_{1})= x_{1},$
 where $1$ denotes the chain
\[\begin{array}{c}
  1\colon\{p\}\to D^{d}_1\times SO(d) \\
 \ \  p\mapsto (T_{0,1},1)\\
\end{array}\]
 Clearly $K_{\theta_{1}(1;x_{1})}=\{p\}
\times e^{-1}(\Delta_{1})=\{p\}\times K_{x_{1}}\cong K_{x_{1}}.$
Moreover
\[ \theta_{1}(1;x_{1})(p;c_{1})(y) =
 \left\{ \begin{array}{ll}
e(x_{1}(c_{1}))  \ \ \mbox{if}\ y
\notin \bigcup \overline{\im(T_{0,1})} \\
          \ &  \\
x_{1}(c_{1}) (y) \ \ \mbox{if}\ y \in \overline{\im(T_{0,1})}
\end{array}  \right.\]
Since $y \in \overline{\im(T_{0,1})}$ for all $y\in D^{d}$ we have
that
$$\theta_{1}[(1;x_{1})(c,c_{1})](y)=x_{1}(c_{1})(y)$$ as it
should, thus axioms 1 holds. Axiom 1 follows from a dimensional
counting argument. Axiom 3 contains two statements, namely, that
the domains and the chain maps associated with both sided of the
commutative diagram agree. The first statement is a consequence of
Theorem \ref{l1}. The second statement follows essentially from
the fact that $M^{S^d}$, the space of smooth from the sphere
sending a neighborhood of the north pole into the fixed point $p
\in M$, is in a natural way a $fD^d$-algebra \cite{BV}.
\\

Next result -- due to Sullivan and Voronov \cite{V1} -- is a
consequence of Theorem \ref{inco}, Theorem \ref{Teo4}  and the
characterization of $H(fD^d)$-algebras given by Salvatore and Wahl
in \cite{sw}. Actually we use the reformulation of the
Salvatore-Wahl theorem given in \cite{CoVo}.

\begin{cor}{\em The graded vector space $\Ho(M^{S^{d+1}})$ is a $H(fD^{d+1})$-algebra, i.e. it is provided
with the following algebraic structures. Let $x,y,z$ be
homogeneous elements of $\Ho(M^{S^d})$.
\begin{enumerate}
\item An associative graded commutative product.

\item A degree $n$ bracket such that $\Ho(M^{S^{d+1}})[n]$ is a graded Lie algebra.

\item $[x, yz]= [x,y]z + (-1)^{(\bar{x}+d)\bar{y}}y[x,z].$

\item For $d$ odd there are operators
$B_i:\Ho(M^{S^{d+1}}) \longrightarrow \Ho(M^{S^{d+1}})[4i-1]$ for
$i $ in $\{1,...,\frac{d-1}{2}\}.$ There is an operator
$\Delta:\Ho(M^{S^{d+1}})
\longrightarrow \Ho(M^{S^{d+1}})[d]$ called the BV operator such that:
\begin{itemize}
\item $ \Delta^2=0.$

\item $(-1)^{\bar{x}}[x,y]=\Delta(xy) -\Delta(x)y - (-1)^{\bar{x}}x\Delta(y)  .$

\item $\Delta[x,y]=[\Delta(x),y] - (-1)^{\bar{x}}[x, \Delta(y)] .$

\end{itemize}
\item For $d$ even there are operators
$B_i:\Ho(M^{S^{d+1}}) \longrightarrow \Ho(M^{S^{d+1}})[4i-1]$ for
$i$ in $\{1,...,\frac{d}{2}\}.$

\item Either in the even or odd case the operators $B_i$ are such that
\begin{itemize}
\item $B_i^2=0.$

\item $B_i$ is a graded derivation on the graded commutative
algebra $\Ho(M^{S^{d+1}})$.

\item $B_i$ is a graded derivation on the graded Lie algebra $\Ho(M^{S^{d+1}})[n].$

\end{itemize}

\end{enumerate}}
\end{cor}

\section{Transversal $1$-categories}\label{sec5}

\begin{figure}[h]
\centering
\includegraphics[width=0.6\textwidth]{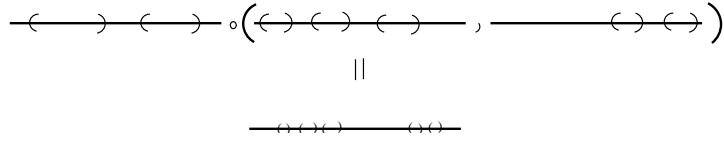} \caption{ \ Composition of little intervals.}\label{Bor9}
\end{figure}

The operad of little discs in dimension $1$ is usually called the
operad of little intervals and is denoted by $I$. Figure
\ref{Bor9} shows an example of composition in the operad of little
intervals. Algebras defined over the operad of little intervals
are called $1$-algebras. It is easy to see that the homology of a
$1$-algebra is an associative algebra. We now introduce the
corresponding notion for the case of pre-categories.

\begin{defi}
 {\em A transversal $1$-category  is a transversal
dg-precategory over the operad $C(I)$ of chains of little
intervals.}
\end{defi}

Let $M$ be a compact manifold and $N_{0}, N_{1}$ be connected
oriented embedded submanifolds of $M.$ Let $Y$ be a smooth compact
manifold.  We denote by $M^{S(Y)}(N_{0},N_{1})$ the set of all
smooth maps $f\colon Y\times [-1,1] \to M,$ such that $f(y,-1)\in
N_{0},$ $f(y,1)\in N_{1},$ and $f$ is a constant map in open
neighborhoods of $Y\times\{-1\}$ and $Y\times\{1\}$, respectively.
$M^{S(Y)}(N_{0},N_{1})$ is a topological space provided with the
compact-open topology. Notice that $M^{S(Y)}(N_{0},N_{1})$ is a
subspace of $\ma(S(Y),M)$ where $$S(Y)= Y\times[-1,1]/\sim$$ and
$\sim$ is the equivalence relation on $Y\times [-1,1]$ given by
$y_{1}\times \{-1\}\sim y_{2}\times\{-1\}$ and $y_{1}\times
\{1\}\sim y_{2}\times\{1\}$ for all $y_{1},y_{2}\in Y.$\\

Consider the complex vector space  $$C(M^{S(Y)}(N_{0},N_{1}))=
\bigoplus_{i=0}^{\infty}
C_{i}(M^{S(Y)}(N_{0},N_{1}))$$  generated by chains $x\colon
K_{x}\longrightarrow M^{S(Y)}(N_{0},N_{1})$ such that the maps
$$e_{-1}(x)\colon K_{x}\times Y\longrightarrow N_{0}
\mbox{ \ \ and \ \ } e_{1}(x)\colon K_{x}\times Y\to N_{1}$$ given by
$e_{i}(x)(c,y)=x(c)(y,i)$ are smooth for $i=-1,1$. Consider the
maps $$e_{-1}\colon M^{S(Y)}(N_{0}, N_{1})\to N_{0} \mbox{ \ \ and
\ \ } e_{1}\colon M^{S(Y)}(N_{0},N_{1})\to N_{1}$$  given respectively
by
$$e_{-1}(f)=f(y,-1)\in N_{0} \mbox{ \ \  and \ \  } e_{1}(f)=f(y,1)\in N_{1}.$$ We
also denote by $e_{i}$ the induced map $e_{i}\colon
\Ca(M^{S(Y)}(N_{0},N_{1})) \longrightarrow \Ca(N_{0})$ given by $$e_{i}(\sum a_{x}x)=\sum
a_{x}e_{i}(x).$$ Given chains $x_{i}\colon K_{x_{i}}\to
M^{S(Y)}(N_{i-1},N_{i})$ for $1 \leq i \leq k,$ consider the map
$$e (x_{1}, \ldots , x_{n})\colon
 \des\prod_{i=1}^{k}K_{x_{i}}
\longrightarrow
\des \prod_{i=1}^{k-1}N_{i}\times N_{i}$$
that sends $(c_{1},\ldots , c_{k})$ to
$$(e_{1}(x_{1}(c_{1})) ,e_{-1}(x_{2}(c_{2})),e_{1}(x_{2}(c_{2})),
\ldots ,e_{-1}(x_{k}(c_{k}))) .$$

Set $$\Omega_{k}=\des\prod_{i=1}^{k}
\Delta^{N_{i}}_{2}\subset\des\prod_{i=1}^{k-1}N{i}\times N{i}$$
where $\Delta^{N_{i}}_{2}=\{(a,a)\in N_{i}\times N_{i} \}$ for
$1\leq i \leq k-1.$ Clearly

$$\!\!\!e^{-1}(\Omega_{k})=\Big\{(c_{1}, \ldots ,c_{k}) \in
\des\prod_{i=1}^{k}K_{x_{i}} \Big |
\begin{array}{c}
e_{1}(x_{i}(c_{i}))=e_{-1}(x_{i+1}(c_{i+1}))
 \\
1\leq i \leq k-1\\
\end{array}\Big \}$$

According to Theorem
\ref{l1} if $e(x_{1},
 \ldots , x_{k})\pitchfork\Omega_{k}$ then
\[e^{-1}(\Omega_{k})=K_{x_{1}}\times_{N_{1}}K_{x_{1}}\times_{N_{2}} \cdots
\times_{N_{k-1}} K_{x_{k}}\] is a manifold with
corners.

\begin{defi}
 {\em  The dg-precategory $\Ca(M^{S(Y)})$ of chains
of branes of type $Y$ in $M$ is such that:
\begin{enumerate}
\item  Objects of $\Ca(M^{S(Y)})$ are connected oriented
embedded submanifolds of $M.$

\item   For $N_{0},N_{1}$ objects in $\Ca(M^{S(Y)})$ we set
$$\Ca(M^{S(Y)})( N_{0},N_{1})=C(M^{S(Y)}( N_{0},N_{1}))[\dim N_1].$$
\end{enumerate}}
\end{defi}

We define a domain in $\Ca(M^{S(Y)})^{\star}$ in $\Ca(M^{S(Y)})$
as follows: given $N_{0},\cdots ,N_{k} \in
\Ca(M^{S(Y)})$ let   $$\Ca(M^{S(Y)}( N_{0},\cdots , N_{k}))\subseteq
\des\bigotimes_{i=1}^{k}\Ca(M^{S(Y)}( N_{i-1},N_{i}))$$
be the space generated by tuples $x_{1}\otimes\cdots
\otimes x_{k}$ such that:

\begin{itemize}
\item $x_{i}\in
\Ca(M^{S(Y)}( N_{i-1},N_{i}))$ for $1\leq i \leq k$,
\item $e(x_{1},\cdots, x_{k})\pitchfork
\Omega_{k}$.
\end{itemize}

The proof that $\Ca(M^{S(Y)})^{*}$ is a domain is similar to the
proof of Lemma \ref{lem 16}.

\begin{thm}\label{t3}
 {\em $\Ca(M^{S(Y)})$  is a
transversal  $1$-category.}
 \end{thm}

Suppose we are given objects $N_{0}, \cdots ,N_{k}$ in
$\Ca(M^{S(Y)})$ we introduce maps
\[\theta_{k}(N_{0}, \cdots ,N_{k})\colon
C(I_k) \otimes \des\bigotimes_{i=1}^{k}
\Ca(M^{S(Y)}(N_{i-1},N_{i}))
\longrightarrow
\Ca(M^{S(Y)})(N_{0},N_{k})\] as follows. Given $x \in C(I_k)$
and $x_{i} \in \Ca(M^{S(Y)}(N_{i-1},N_{i})),$ then
$\theta_{k}(x;x_{1},\ldots ,x_{k})$ has domain
$$K_{\theta_{k}(x;x_{1},\ldots ,x_{k})}=K_{x}
\times e^{-1}(\Omega_{k} ).$$
Let $x\colon K_{x}\longrightarrow I_k$ be such that for $c\in
K_{x}$ we have
$$x(c)=(T_{p_{1}(c),r_{1}(c)},\cdots , T_{p_{k}(c),r_{k}(c)}).$$ The map $$\theta_{k}(x;x_{1},\ldots ,x_{k})\colon
K_{\theta_{k}(x;x_{1},\ldots ,x_{k})} \longrightarrow \Ca(
M^{S(Y)})$$ is such that for $t\in I$ and $y \in Y$ we have
\[\theta [{(x;x_{1},\ldots ,x_{k})(c;c_{1},\ldots ,c_{n})](y,t)}
 = \left\{ \begin{array}{ll}
e_{1}(x_{i}(c_{i}))  \ \ \mbox{if} \ t \notin \bigcup
\overline{\im(T_{p_{i}(c),r_{i}(c)})}\\
 \ & \ \\
x_{i}(c_{i})\Big(y,\frac{t-p_{i}(c)}{r_{i}(c)}\Big) \ \ \mbox{if}\
t \in \overline{\im(T_{p_{i}(c),r_{i}(c)})}

          \ & \ \\
\end{array} \right. \]
Axioms 1,2,3 of Definition \ref{d2} are  proved as Theorem
\ref{inco}.\\

\begin{figure}[h]
\centering
\includegraphics[width=0.5\textwidth]{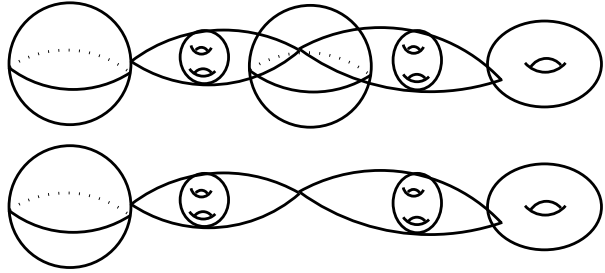} \caption{ \ Schematic representation of compositions
in  $\Ca(M^{S(Y)})$. }\label{Bor99}
\end{figure}

Figure \ref{Bor99} represents, schematically, the composition in
the category $\Ca(M^{S(Y)})$ where $Y$ is a surface of genus
$2$.\\

\begin{figure}[h]
\centering
\includegraphics[width=0.5\textwidth]{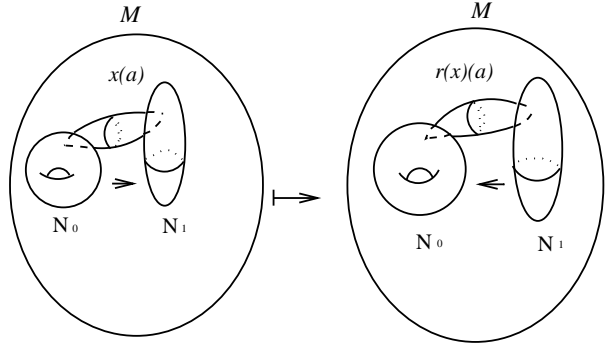} \caption{ \ Example of an application of the functor $\ri$.}\label{Bor7}
\end{figure}

Let us consider the map $^{-}\colon I_n\longrightarrow I_n$ given
by $$\overline{(T_{a_{1}r_{1}}, \cdots , T_{a_{n}r_{n}})}=
(T_{-a_{1}r_{1}}, \cdots ,T_{-a_{n}r_{n}}).$$ We also need the
induced chain map $^{-}\colon
\Ca(I(n))\longrightarrow
\Ca(I(n)).$ An interesting feature of the $1$-category
 $\Ca(M^{S(Y)})$ is that it comes with a natural contravariant
prefunctor  $$\ri\colon \Ca(M^{S(Y)})\longrightarrow
\Ca(M^{S(Y)})$$  which is the identity on
 objects; for objects $N_{0},N_{1}$ in  $\Ca(M^{S(Y)})$ the map
\[\ri \colon \Ca(M^{S(Y)})(N_{0},N_{1})\longrightarrow
\Ca(M^{S(Y)})(N_{1},N_{0})\] is defined  as follows: for $x\in
\Ca(M^{S(Y)})(N_{0},N_{1})$
the domain of $\ri(x)$  is $K_{x}$ and if $c\in K_{x}$ then for
$-1\leq t \leq 1$ we set
$$[\ri(x)(c)](y,t)=[x(c)](y,-t).$$  Figure \ref{Bor7} illustrates the
meaning of the functor $\ri.$ It is not hard to check that $\ri$
satisfies the following identity
\[\ri (\theta_{n}(x;x_{1},\cdots ,x_{n}))= \pm
 \theta_{n}(\overline{x};
 \ri( x_{n}),\cdots ,\ri(x_{1}))\circ s\]
where $x\in C(I_n)$, $x_{i}\in
 \Ca(M^{S(Y)})(N_{i-1},N_{i})$  and the map
 \[s\colon K_{x}\times K_{x_{1}}\times_{N_{1}}K_{x_{2}}\times \cdots  \times_{N_{n-1}}
K_{x_{n}} \longrightarrow  K_{x}\times
K_{x_{n}}\times_{N_{n-1}}K_{x_{n-i}}\times \cdots  \times_{N_{1}}
K_{x_{1}}\] is given for $a\in K_{x}$ and $ t_{i}\in K_{x_{i}}$ by
$$s(a;t_{1},\cdots t_{n})=(a;t_{n},\cdots t_{1}).$$

We need some notions from universal algebra. The concepts that we
need where introduce by Markl in
\cite{MM}, where the reader will find further details.

\begin{figure}[h]
\centering
\includegraphics[width=0.4\textwidth]{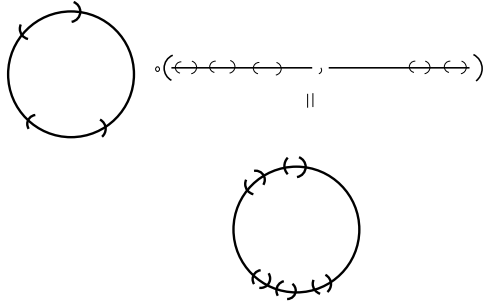} \caption{ \ Example of the right
$I_n$-module structure on $S^{1}_n$. }\label{Bor11}
\end{figure}

\begin{defi}\label{Tramod}
 {\em A   right $C(I)$-module $M$
consists of a sequence  $M_n$  of objects in dg-vect together with
maps for $k\geq 0$ $$\lambda_{k}\colon M_k\otimes \bigotimes^{
k}_{s=1} C(I_{j_{s}}) \longrightarrow M{j_{1}+\cdots +j_{k}}$$
that are associative and unital.}
\end{defi}
Consider the space $S^{1}_n$ of configurations of $n$ little discs
inside the unit circle. $S^{1}_n$ is obtained from $I_n$ by
identifying the ends points of the interval $[-1,1].$ Markl in
\cite{MM} shows that $S^{1}_n$ is a right $I_n$-module in the
topological category, as usual that result implies the following
result. The compositions given $S^{1}_n$ the structure of a right
$I_n$-module is illustrated in Figure \ref{Bor11}.

\begin{lema}\label{lem29}
 {\em $C(S^{1}_n)$ is a right $C(I_n)$-module.}
\end{lema}

\begin{defi}\label{tra}
 {\em Let $\C$ be a transversal $1$-category. A  $C(S^{1})$-trace over $\C$
is an object $B$ in dg-vect together with maps
\[ T_{N_{0}, \cdots ,N_{k-1}}\colon C(S^{1}_k) \otimes \C(N_{0},\cdots
,N_{k-1},N_{0})\longrightarrow B\] for $N_{0}, \cdots ,N_{k-1}$
objects of $\C,$ such that the following diagram is commutative}
\[\xymatrix
 { C(S^{1}_k) \otimes \des\bigotimes_{i=1}^{k}
C(I_{j_{i}}) \otimes \C(N_{0},\cdots ,N_{j-1},N_{0})
 \ar[r]^{\ \ \ \ \ \ \ \ \ \ \ \ \lambda_{k}\otimes\id}
\ar @ {^{(}->}[d]  & \Ca(S^{1}_j)\otimes
\C(N_{0},\cdots N_{j-1},N_{0})
\ar[d]^{T }\\
C(S^{1}_k) \otimes
\des\bigotimes_{i=1}^{k}C(I_{j_{i}})\otimes
\des\bigotimes_{i=1}^{k} \C(N_{j_{i-1}},\cdots
N_{j_{i}})\ar[d]_{\mbox{shuffle}} & B
  \\
C(S^{1}_k) \otimes
\des\bigotimes_{i=1}^{k}C(I_{j_{i}})\otimes\des\bigotimes_{i=1}^{k}
\C(N_{j_{i-1}},\cdots N_{j_{i}})
 \ar[r] &
C(S^{1}_k)
\otimes \des\bigotimes_{i=1}^{k}\C(N_{j_{i-1}},N_{j_{i}})
\ar[u]_{T} }\]
\end{defi}

Let $M$ be a compact oriented smooth manifold and  $Y$ be a
compact smooth manifold. We denote by $M^{Y\times S^{1}}$ the set
of  smooth maps $f\colon Y\times S^{1} \to M.$  We impose on
$M^{Y\times S^{1}}$  the compact-open topology and we set
$$C(M^{Y\times S^{1}})= \des\bigoplus_{i=0}^{\infty}
C_{i}(M^{Y\times S^{1}})$$  where $\Ca_{i}(M^{Y\times S^{1}})$ is
the vector space generated by chains $x\colon K_{x}\longrightarrow
M^{Y\times S^{1}}$ such that the induced map
$$\widehat{x}\colon K_{x}\times Y\times S^{1}\longrightarrow M$$ is a smooth map.

\begin{thm}\label{to3}
 {\em $\Ca(M^{S(Y)})$ admits a natural $C(S^{1})$-trace.}
\end{thm}

To prove this result we define maps
$$T\colon
 C(S^{1}_n)\otimes \Ca(M^{S(Y)})(N_{0},N_{1}) \otimes \cdots \otimes
\Ca(M^{S(Y)})(N_{k-1},N_{0}) \longrightarrow \Ca(M^{Y\times S^{1}}).$$
Assume we are given chains $$x \in C(S^{1}_k), \ \ x_{i} \in
\Ca(M^{S(Y)}(N_{i-1},N_{i})) \mbox{ \ and \ } x_{k}\in
\Ca(M^{S(Y)}(N_{k-1},N_{0})).$$ The $k$-tuple $x_{1}\otimes\cdots
\otimes x_{k}$ belongs to  $\Ca(M^{S(Y)})( N_{0},\cdots
,N_{k-1},N_{0})$ if and only if $$e(x_{1},\cdots x_{k})\pitchfork
\Omega_{k}.$$
Then $T(x;x_{1},\ldots ,x_{k})$ is the chain with domain
\[K(x;x_{1},\ldots
,x_{k})=K_{x}\times e^{-1}(\Omega_{k})\] and $$T(x;x_{1},\ldots
,x_{k}) \colon K_{x}\times
 K_{x_{1}}\times_{N_{1}}K_{x_{1}}\times_{N_{2}} \cdots \times_{N_{0}}
 K_{x_{0}} \longrightarrow  M^{Y\times S^{1}}$$  is the map given by
\[ [T(x;x_{1},\ldots ,x_{k})(c,c_{1},\ldots ,c_{k})](t) = \left\{ \begin{array}{ll}
e_{1}(x_{i}(c_{i}))  \ \ \mbox{if} \ t \notin \bigcup
\overline{\im(T_{p_{i}(c),r_{i}(c)})}\\
 \ & \ \\
x_{i}(c_{i})\Big(\frac{t-y_{i}(c) }{r_{i}(c)}\Big)
\ \ \mbox{if}\ t \in \overline{\im(T_{y_{i}(c),r_{i}(c)})}  \\
\end{array} \right. \]
The rest of the proof is similar to that of Theorem \ref{inco}.\\

\begin{lema}\label{to7}
 {\em The category $\Ho(M^{S(Y)})$ admits a natural
$H(S^{1})$-trace.}
\end{lema}

The result follows from Lemma \ref{lem29}, Theorem \ref{Teo4} and
Corollary
\ref{cor5}.\\

Note that in the case that $Y$ is a point we recover known results
from open string topology \cite{DS}. The category of homological
open string carries additional structures, for example Sullivan
has defined a co-category structure on it, and more generally
Baas, Cohen and Ramirez \cite{baas} have shown that there are
further categorical operations coming from surfaces of higher
genera with boundaries and  marked intervals on them. Tamanoi has
discusses in
\cite{tama1} the conditions on the surfaces such that the
corresponding operations are not necessarily trivial.  It is not
clear to us if these additional structures are also present on
$\Ho(M^{S(Y)})$ with $Y$ a positive dimensional manifold.

\section{Homological Quantum Field Theory}\label{sec7}

In this section we shall introduce the main definition of this
work, namely, the notion of homological quantum field theories. To
understand this notion two prerequisites are needed: the string
topology of Chas and Sullivan that we have discussed in the
previous sections, and the categorical approach \cite{MA} towards
quantum field theory which we proceed to review. \\

\noindent \textbf{Category $\Co_d$ of $d$ dimensional cobordisms.} The leading role
in the categorical approach to quantum field theory is the
category  $\Co_d$ of cobordism introduced by Ren\'e Thom in
\cite{RT, RT1}.
For $d \geq 1$ objects in $\Co_d$ are  compact oriented $d-1$
dimensional smooth manifolds. Morphisms between objects $N_1$ and
$N_2$ in $\Co_d$ are of two types:\\

\noindent 1) A diffeomorphism from $N_1$ to $N_2$.\\

\noindent  2) A compact oriented
manifold with boundaries $M$ -- a cobordism -- together with a
diffeomorphism from $ (-N_1 \sqcup N_2) \times [0,1)$ onto an open
neighborhood of $\partial(M) $. This kind of morphisms are
considered up to diffeomorphisms.\\

Composition of morphisms $\Co_d$ is given by composition of
diffeomorphisms and gluing of cobordisms manifolds. It is usually
assumed that both morphisms and objects are provided with
additional data. Thus we postulate that there is in addition a
contravariant functor $$D: OBMan_d \longrightarrow Set$$ from the
category whose objects are either $d$ dimensional manifolds with
boundaries, or $d-1$ manifolds without boundaries. Morphisms in
$OBMan_d$ are smooth maps. With the help of the functor $D$ we
define a colored category of cobordisms $D\Co_d$ whose objects are
pairs $(N,s)$ where $N$ is an object of $\Co_d$ and $s \in F(N).$
We think of $s$ as giving additional structure to the manifold
$N$. Morphisms of type 1) are structure preserving
diffeomorphisms. Morphisms of type 2) are  cobordisms pairs
$(M,s)$ where $M$ is a cobordism and $s \in F(M)$. It is required
that the structure $s$ when restricted to the boundary of $M$
agrees with the structure originally given to the boundary
components of $M$.  \\

\noindent \textbf{Monoidal representations of $\Co_d$}. Gradually it has become clear that
the geometric background  for the mathematical understanding of
quantum fields is given by monoidal representations of the
category of structured cobordisms, i.e., monoidal functors
$$F:D\Co_d
\longrightarrow \ve$$ from $D$-cobordisms into vector spaces. Field theories are not
determined by its geometric background and there are additional
constrains for a realistic quantum field theory than those imposed
by the fact that they yield monoidal representations of $D\Co_d $.
Different types of field theories correspond to different  choices
of different types of data on the objects and morphisms of the
cobordisms category, i.e. different choices of the functor $D$. It
is often the case that the sets of morphisms in $D\Co_d$ come with
a natural topology. In those cases, a field theory is a continuous
monoidal functor from $D\Co_d$ into $\ve$. Some of the most
relevant types of theories from this point of view are the
following:

\begin{itemize}
\item \textit{Lorentzian quantum field theory} $LQFT$. For a rather comprehensive mathematical introduction
to field theory the reader may consult \cite{de}. Unfortunately,
the analytical difficulties have prevented, so far, fully rigorous
constructions of field theories of this type. One considers the
category $L\Co_d$ of Lorentzian cobordisms defined as $\Co_d$ with
the extra data: objects are provided with a Riemannian metric,
morphisms are provided with a Lorenzian metric such that its
restriction to the boundary components agree with the specified
Riemannian metric on objects. Lorentzian quantum field theories
$LQFT$ are linear representation of the category $L\Co_d$, i.e.,
monoidal functors $F:L\Co_d
\longrightarrow \ve$.

\item \textit{Euclidean quantum field theory} $EQFT$. One constructs the
category  $E\Co_d$ of Euclidean or Riemannian  cobordisms as in
the Lorentzian situation; in this case the metrics on both objects
and morphisms are assumed to be Riemannian. A Euclidean quantum
field theory $EQFT$ is a monoidal functor $F:E\Co_d
\longrightarrow \ve$.

\item \textit{Conformal field theory} $CFT$.   Riemaniann metrics $g$
and $h$ on a manifold $M$ are said to be conformally equivalent if
there exist a diffeomorphism $f:M \longrightarrow M$ and a smooth
map $\lambda:M \longrightarrow \mathbb{R}_+$ such that
$f^{*}(g)=\lambda h.$ The category of conformal cobordisms
$C\Co_d$  is defined as in the Euclidean case but now we demand
that  objects and morphisms be provided with Riemaniann metrics
defined up to conformal equivalence. Conformal field theories
$CFT$ are monoidal functors $F: C\Co_d \longrightarrow \ve$.
Unlike the previous types this sort of theory has been deeply
studied in the mathematical literature. Kontsevich in \cite{KO1}
has proposed that conformal field theories are deeply related with
$d$-algebras. The case $d=2$ was first axiomatized by Segal in
\cite{Se}, it has attracted a lot of attention because of its
relation with string theory, and because this case may be treated
with complex analytic methods since a conformal metric on a
surface is the same as a complex structure on it. There have been
many developments in the subject out of which we cite just a few
\cite{and, cos1, cos2, DP0}.

\item \textit{Topological quantum field theory} TQFT.  This sort of
theory was described within the framework of linear
representations of $\Co_d$ by Atiyah in \cite{MA, MA2}. In a sense
this sort of theory is the prototype that indicates the
possibilities of the categorical approach; it has been deeply
studied in the literature, for example in the works
\cite{ale, Baez, dij, kim1, law, Ur, pick, q, Toen, TU, wk}. In essence,  the
category  TQFT of topological quantum field theories may be
identified with the category $\Mfu(\Co_{n},\ve)$ of monoidal
functors $$F \colon
\Co_{n}\longrightarrow
\ve,$$ that is, topological quantum field theory deals with the
bare category of cobordisms without further structures imposed on
its objects or morphisms.

\item \textit{Homotopical quantum field theory} HQFT. This sort of
theory was introduced by Turaev in \cite{TU, TU3} and has been
further developed, among others, by Brightwell, Bunke, Porter,
Rodrigues, Turaev, Turner and Willerton  \cite{bri, bun, por1,
port2, rod}. Fix a compact connected smooth manifold $M$. The
category  $H\Co_d^{M}$ of homotopically extended cobordisms $M$ is
such that its objects are $d-1$ dimensional smooth compact
manifolds $N$ together with a homotopy class of maps $f:N
\longrightarrow M$. Morphisms in $H\Co_d^{M}$ from $N_0$ to $N_1$ are
cobordisms $P$ connecting $N_0$ and $N_1$ together with a homotopy
class of  maps $g:P
\longrightarrow M$ such that its restriction to the boundaries gives
the agrees with the homotopy classes associate with them. A
homotopical quantum field theory is a monoidal functor
$F:H\Co_d^{M} \longrightarrow \ve$.\\

\item  \textit{Homological quamtun field theory} HLQFT.  A complete
definition of this sort of theory is the main topic of this
section. First one construct the category $\Co_d^{M}$ of
homological extended cobordisms. Its objects are $d-1$ dimensional
manifolds $N$ together with a map sending each boundary component
of $N$ into an oriented embedded submanifolds of $M$. Morphisms
are cobordisms together with an homology class of maps (constant
on a neighborhood of each boundary component and mapping each
boundary component into its associated embedded submanifold)  from
the cobordism into $M$. The compositions in $\Co_d^{M}$ are
defined in a rather interesting fashion using techniques
originally
introduced Chas and Sullivan  in the context of String topology.\\

\end{itemize}

We proceed to define in details the category $\HL_d$. It would be
done in the following steps:

\begin{itemize}

\item We construct a transversal $1$-category $\sC \OB_{d}^{M}$ for each integer $d\geq1$ and
each compact oriented smooth manifold $M.$

\item $\Co_{d}^{M}$ is defined by the identity
$\Co_{d}^{M}=\Ho(\sC\OB_{d}^{M}).$

\item  $\HL_d(M)$ is defined as the category   $\Mfu(\Co_{n}^M,\ve)$ of monoidal
functors from $\Co_{d}^{M}$ to $\ve.$
\end{itemize}

Objects of $\sC \OB_{n}^{M}$ are triples $(N,f,<)$ such that:
\begin{enumerate}
\item $N$ is a compact oriented  manifold of dimension $d-1.$

\item $f\colon \pi_{0}(N) \longrightarrow D(M)$  is any map. For
such a map $f$ we set $\overline{f}=\prod_{c \in
\pi_{0}(N)}f(c).$

\item $<$ is a linear ordering on $\pi_{0}(N).$
\end{enumerate}

By convention the empty set is assumed to be a $d$-dimensional
manifold for all $d\in
\mathbb{N}.$ Let $(N_{0},f_{0},<_{0})$ and $(N_{1}, \ f_{1},<_{1})$ be objects in
$\sC \OB_{d}^{M}$, we set
\[ \sC \OB_{d}^{M}((N_{0},f_{0},<_{0}), \ (N_{1},f_{1},<_{1}))=\overline{\sC
\OB_{d}^{M}}\diagup \backsim ,\] where by definition $\overline{\sC \OB_{d}^{M}}$
is the set of triples $(P,\alpha,\xi)$ such that:
\begin{itemize}

\item $P$ is a compact oriented smooth manifold with corners of
dimension $d.$

\item $\alpha \colon N_{0}\bigsqcup N_{1}\times
[0,1]\longrightarrow
\im (\alpha)\subseteq P$ is a diffeomorphism and $\alpha \mid_{
N_{0}\bigsqcup N_{1}}\longrightarrow
\partial  P$ is such that $\alpha|_{N_{0}}$
reverses orientation, and $\alpha|_{N_{1}}$ preserves orientation.

\item $\xi \in \Ca(M^P_{f_{0}, f_{1}})=C(M^P_{f_{0}, f_{1}})[\dim N_{1}],$ where $M^P_{f_{0},
f_{1}}$ is the space of smooth maps $g\colon P\longrightarrow M$
such that for each $c\in
\pi_{0}(N_{j}),$ $g$ is a constant map with value in
$f_{j}(c)$ on an open neighborhood of $c,$ for $j=0,1.$
\end{itemize}
For $g \in M^P_{f_{0},f_{1}}$ we define $e_{0}(g)\in
\overline{f}_{0}$,  $e_{1}(g)\in
\overline{f}_{1}$ by $e_{0}(g)(c)=e_{0}(x)$, $e_{1}(g)(c)=e_{1}(x)$ for any $x \in c.$ We define an
equivalence relation on $$\overline{\sC
\OB_{d}^{M}}((N_{0},f_{0},<_{0}), \ (N_{1},f_{1},<_{1}))$$ as follows: triples $(P_{1},\alpha_{1},
\xi_{1} )$ and $(P_{2},\alpha_{2}, \xi_{2} )$ are equivalent if
there is an orientation preserving diffeomorphism $\varphi\colon
P_{1} \longrightarrow P_{2}$ such that $\varphi \circ
\alpha_{1}=\alpha_{2},$ and $\varphi_{\star}(\xi_{1})=\xi_{2}.$\\

Let $\sC \OB_{d}^{M} ((N_{0},f_{0},<_{0}),\cdots ,
(N_{k},f_{k},<_{k}))$ be the vector space generated by $k$-tuples
$\{(P_{i},\alpha_{i},
\xi_{i})\}_{i=1}^{k}$ such that for $1 \leq i \leq k$ we have that
$$(P_{i},\alpha_{i},\xi_{i})\in \sC \OB_{d}^{M}
((N_{i-1},f_{i-1},<_{i-1}), (N_{i},f_{i},<_{i})),$$ and the map
$$e (\xi_{1},
 \ldots , \xi_{k})\colon
 \des\prod_{i=1}^{k}K_{\xi_{i}}
\longrightarrow
\des \prod_{i=1}^{k-1}\overline{f_{i}}\times\overline{ f_{i}}$$
given by $$e (\xi_{1}, \ldots , \xi_{k}) (c_{1},\ldots , c_{k}) =
(e_{1}(\xi_{1}(c_{1}))
,e_{0}(\xi_{2}(c_{2})),e_{1}(\xi_{2}(c_{2})),
\ldots ,e_{0}(\xi_{k}(c_{k})))$$ is transversal to $\Omega_{k}=\prod_{i=1}^{k-1}
\Delta^{\overline{f_{i}}}_{2}\subset \prod_{i=1}^{k-1}\overline{f_{i}}\times
\overline{f_{i}}$ where for $1\leq i \leq k-1$ we set
$$\Delta^{\overline{N_{i}}}_{2}=\{(a,a)\in
\overline{f_{i}}\times \overline{f_{i}} \}.$$
Clearly we have that $$ e^{-1}(\Omega_{k})= \Big\{(c_{1}, \ldots
,c_{k}) \in
\des\prod_{i=1}^{k}K_{\xi_{i}} \Big |
\begin{array}{c}
e_{1}(\xi_{i}(c_{i}))=e_{0}(\xi_{i+1}(c_{i+1}))
 \\ 1\leq i \leq k-1\\
\end{array}\Big \}$$
Since $e$ is a smooth map and $e(\xi_{1},
 \ldots , \xi_{k})\pitchfork\Omega_{k}$ then
\[e^{-1}(\Omega_{k})=K_{\xi_{1}}\times_{\overline{f}_{1}}K_{\xi_{1}}\times_{\overline{f}_{2}} \cdots
\times_{\overline{f}_{k-1}} K_{x_{k}}\] is a
manifold with corners.\\

\begin{figure}[h]
\centering
\includegraphics[width=0.4\textwidth]{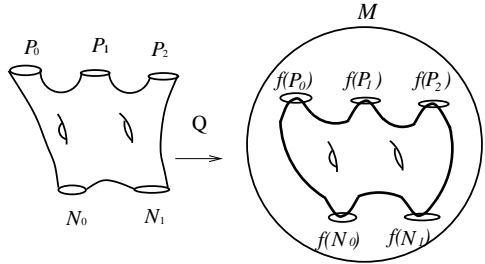} \caption{ \
Example of a $d$-cobordism enriched over $M.$}\label{Bor38}
\end{figure}

\begin{figure}[h]
\centering
\includegraphics[width=0.4\textwidth]{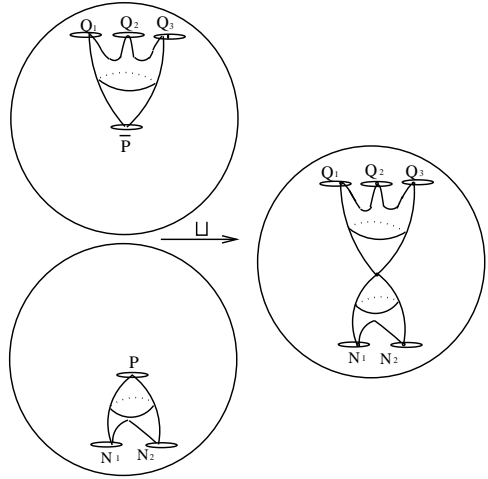} \caption{ \ Composition of
$d$-cobordism enriched over $M.$}\label{Bor49}
\end{figure}

Given $a\in C(\I_k)$ and chains $(P_{i},\alpha_{i},\xi_{i})\in
 \sC \OB_{d}^{M}
((N_{i-1},f_{i-1},<_{i-1}), (N_{i}, f_{i}, <_{i}))$ for $1\leq
i\leq k,$ the composition morphism
$$a((P_{1},\alpha_{1},\xi_{1}),\cdots
,(P_{k},\alpha_{k},\xi_{k}))\in \sC \OB_{n}^{M} ((N_{0},
f_{0},<_{0}), (N_{k},f_{k},<_{k}))$$ is the triple
$(a(P_{1},\cdots , P_{k}),a(\alpha_{1},\cdots ,
\alpha_{k}),a(\xi_{1},\cdots , \xi_{k}))$ such that
\[\begin{array}{c}
  a(P_{1},\cdots , P_{k})=P_{1}\des\bigsqcup_{N_{1}}
\cdots \des\bigsqcup_{N_{k-1}}P_{k} \\
  a(\alpha_{1},\cdots ,
\alpha_{k})=\alpha_{1}\mid_{N_{0}}\bigsqcup
\alpha_{k}\mid_{N_{k}} \\
          \\
  K_{a}(\xi_{1},\cdots , \xi_{k})=K_{a}\times
  K_{\xi_{1}}\times_{\overline{f}_{1}}K_{\xi_{2}}\times
\cdots \times_{\overline{f}_{k-1}}K_{\xi_{k}} \\
\end{array}\]
The map $a(\xi_{1},\cdots , \xi_{k})\colon K_{a}(\xi_{1},\cdots ,
\xi_{k})\times P_{1}\sqcup_{N_{1}}\cdots
\sqcup_{N_{k-1}}P_{k}\longrightarrow M$ is given by $$a(\xi_{1},\cdots ,
\xi_{k})(s, t_{1}\cdots ,t_{k},u)=\xi_{i}(t_{i})(u)$$ for
$(s, t_{1}\cdots ,t_{k})$ in $K_{a}\times
K_{\xi_{1}}\times_{\overline{f}_{1}} \cdots
\times_{\overline{f}_{k-1}}  K_{\xi_{k}}$ and $u\in P_{i}.$
Figure $\ref{Bor38}$ represents a $d$-cobordism enriched over $M$
and Figure $\ref{Bor49}$ shows a composition of $d$-cobordism
enriched over $M.$\\

Let $\Co^{M}_{d,r}$  be the full subcategory of $\Co^{M}_{d}$ such
that the empty set is not longer accepted as a valid   object.

\begin{prop}
 {\em $\Co_{d}^{M}$ is a monoidal category
with disjoint union $\sqcup$ as product and empty set as unit.
Furthermore $\Co^{M}_{d,r}$
 is a monoidal category without unit. }
\end{prop}

Given monoidal categories $\C$ and $\d$ we let $\Mfu(\C,\d)$ be
the category of monoidal functors from $\C$ to $\d.$
\begin{defi}
{\em The category of homological quantum field theories of
dimension $d$ is $$\HL_d (M)=\Mfu(\Co_{d}^{M},\ve).$$ The category
of restricted homological quantum field theories of dimension $d$
is $$\HL_{d,r}(M)=\Mfu(\Co_{d,r}^{M},\ve).$$}
\end{defi}

Let us try to digest the meaning of the previous definition.
 Objects in $\HL_d(M)$ are monoidal
functors $F \colon \Co_{d}^{M} \longrightarrow \ve,$ i.e. $F$
assigns to each triple $(N,f,<)$ a vector space $F(N,f,<)$ in such
a way that
$$F(N,f,<)=\bigotimes_{c \in
\pi_0(N)}F(c, f(c), <) .$$ $F$ assigns to each homological cobordism
$\alpha$ from $N$ to $L$ a linear map
$$F(\alpha): \bigotimes_{c \in \pi_0(N)}F(c, f(c), <) \longrightarrow \bigotimes_{c \in \pi_0(L)}F(c, f(c), <).$$
 Morphisms in $\HL_d(M)(F ,G)$ are natural transformations $T\colon F \longrightarrow
G,$ i.e. for each triple $(N,f,<)$ there is a linear map $$T_N:
\bigotimes_{c \in
\pi_0(N)}F(c, f(c), <) \longrightarrow \bigotimes_{c \in
\pi_0(N)}G(c, f(c), <),$$ such that if $\alpha$ is an homologically
extended cobordism from $N$ to $L$ then the following diagram is
commutative

\[\xymatrix @R=.4in  @C=.6in
{\bigotimes_{c \in \pi_0(N)}F(c, f(c), <) \ar[d]^{F(\alpha)}
\ar[r]^{T_N} & \bigotimes_{c \in \pi_0(N)}G(c, f(c), <)\ar[d]^{G(\alpha)}\\
\bigotimes_{c \in \pi_0(L)}F(c, f(c), <) \ar[r]^{T_L} & \bigotimes_{c \in \pi_0(L)}G(c, f(c), <)}\]

There is a canonical restricted homological quantum field theory
attached to any given manifold $M$. Consider the prefunctor $\Ho
\colon \Co_{d,r}^{M}
\longrightarrow
\mbox{vect}$ given on objects by
\[\xymatrix @C=.1in  @R=.02in
{\Ho \colon \ob(\Co_{d}^{M})
 \ar[rr] & &
\ob(\mbox{vect})\\
(N,f,<) \ar @{|->}[rr] & & \
 \Ho(N,f,<)= \Ho(\overline{f})
 }\]
The image under $\Ho$ of $(P,\alpha, \xi)\in
\Co_{d,r}^{M}((N_{0},f_{0},<_{0}), (N_{1}, \ f_{1},<_{1}))$  is
the linear map
\[\xymatrix @C=.1in  @R=.02in {\Ho(P,\alpha, \xi)
    \colon \Ho(N_{0},f_{0},<_{0}) \ar[rr] & &
\Ho(N_{1},f_{1},<_{1})\\
\ \ \ \ \ \ \ \ \ \ \ \ \ \ \ x \ar @{|->}[r] & &
 \Ho(P, \alpha,\xi)(x)}\]
where for $\xi \in \Ho(\ma(P,M)_{f_{0},f_{1}})$ given by
$\xi\colon K_{\xi}\longrightarrow \ma(P,M)_{f_{0}, f_{1}},$ and
$x\in
\Ho(\overline{f}_{0})$ given by $x\colon K_{x}\longrightarrow
\overline{f}_{0},$ the domain $K_{\Ho(P,\alpha ,
\xi)(x)}$ of $\Ho(P, \alpha ,\xi)(x)$ is given by
$$K_{\Ho(P,\alpha,\xi)(x)}=K_{x}\times_{\overline{f}_{0}}
K_{\xi},$$ and  the map $\Ho(P, \alpha ,\xi)(x)$ is given by
$$\Ho(P, \alpha ,\xi)(x)(a,t)=e_{1}(t),$$ where
$t=(t_{c})_{c\in\pi_{0}(f_{0})}.$

\begin{thm}
 {\em $\Ho$ is a restricted homological quantum field theory.}
\end{thm}

Indeed let $$(P_{1},\alpha_{1}, \xi_{1})\in
\Co_{d,r}^{M}((N_{0},f_{0},<_{0}), (N_{1}, \ f_{1},<_{1}))$$ and
$$(P_{2},\alpha_{2}, \xi_{2})\in \Co_{n,r}^{M}((N_{1},f_{1},<_{1}),
(N_{2}, \ f_{2},<_{2})).$$ We must check that
\[\Ho((P_{2},\alpha_{2}, \xi_{2})\circ (P_{1},\alpha_{1},
\xi_{1}))=\Ho(P_{2},\alpha_{2}, \xi_{2})\circ\Ho(P_{1},\alpha_{1},
\xi_{1}).\] Since the domain of $(P_{2},\alpha_{2}, \xi_{2})\circ
(P_{1},\alpha_{1}, \xi_{1})$ is
$K_{\xi_{1}}\times_{\overline{f}_{1}}K_{\xi_{2}},$ then the domain
of
$$\Ho((P_{2},\alpha_{2}, \xi_{2})\circ (P_{1},\alpha_{1},
\xi_{1}))(x)$$ is given by
$$K_{x}\times_{\overline{f}_{0}}(K_{\xi_{1}}\times_{\overline{f}_{1}}K_{\xi_{2}}).$$
 On the other hand the domain of
$\Ho(P_{2},\alpha_{2}, \xi_{2})\circ \Ho(P_{1},\alpha_{1},
\xi_{1})(x)$ is
$$(K_{x}\times_{\overline{f}_{0}}K_{\xi_{1}})\times_{\overline{f}_{1}}K_{\xi_{2}}.$$
Thus we see that the domains agree and it is easy to check that
the corresponding functions also agree.

\section{1-dimensional homological quantum field theories}\label{sss}

In this section, based on \cite{cd3}, we study  examples of
restricted homological quantum field theories in dimension one.
First we show that there is a intimate relationship between
$\Co_{1,r}^{M}$ the category of homologically extended
$1$-dimensional cobordisms and the category $\Ho(M^{I})$ of open
strings \cite{DS} in $M$. This relationship should not be confused
with the fact, due to Cohen-Godin
\cite{Ch3}, that string homology is a restricted two dimensional topological
quantum field theory. Second we show that there are plenty of
non-trivial examples of \HL \ in dimension one, indeed  we show
that one can associate such an object to each connection on
principal fiber bundle. Third we explore the notion of homological
matrices and discuss its relationship with homological quantum
fields theories in dimension one.\\

Let us first show how $1$-dimensional homological quantum field
theories are related to open string topology  which was considered
in Section \ref{sec5} in the case that $Y$ is a point. Objects in
the open string category $\Ho(M^I)$  are embedded submanifolds of
$M$. The space of morphisms from $N_0$ and $N_1$ is given by
$$\Ho(M^{I}_{N_0,N_1})=H(M_{N_0, N_1}^{I})[\dim(N_1)]$$ the
homology with degrees shifted down by $\dim N_1$ of the space
$M_{N_0, N_1}^{I}$ of smooth path $$x:I
\longrightarrow M$$  constant on  neighborhoods of
$0$ and $1$. Composition of morphisms defined in Section
\ref{sec5} yields a map  $$
\Ho(M_{N_0, N_1}^{I})\otimes \Ho(M_{N_1, N_2}^{I}) \longrightarrow
\Ho(M_{N_0, N_2}^{I}) $$ turning $\Ho(M^{I})$ into a graded
category.\\

Let us now consider the category $\HL_{1,r}$ of restricted
homological quantum field theories in dimension $1$, which is
given by
$$\HL_{1,r}=\Mfu(\Co^{M}_{1,r},
\ve) .$$  An object $f$ in $\Co^{M}_{1,r}$ is just a map
$f
\colon [n]\longrightarrow D(M)$ where we set $[n]=\{1,\cdots ,n\}.$  Let $S_n$ be the group of permutations of $n$
letters. We shall use the notation $\overline{f}=
\prod_{i\in[n]}f(i).$  The space of morphisms in $\Co_{1,r}^M$
from $f$ to $g$ is by definition given by
\[\Co^{M}_{1,r}(f,g)= \bigoplus_{\sigma\in S_n} \bigotimes_{i=1}^{n}
\Ho(M^I_{f(i),g(\sigma(i))}).\] Composition of morphisms in $\Co_{1,r}^{M}$ is given by the following composition of
maps:

\[\xymatrix @R=.15in{\Co^{M}_{1,r}(f,g)\otimes \Co^{M}_{1,r}(g,h) \ar[d]\\
\bigoplus_{\sigma, \tau\in S_n} \bigotimes_{i=1}^{n}
\Ho(M^I_{f(i),g(\sigma(i))})\otimes \bigotimes_{j=1}^{n}
\Ho(M^I_{g(j),h(\tau(j))}) \ar[d]\\
 \bigoplus_{\sigma, \tau\in S_n} \bigotimes_{i=1}^{n}
\Ho(M^I_{f(i),g(\sigma(i))})\otimes
\Ho(M^I_{g(\sigma(i))}, M^I_{h(\tau(\sigma(i)))})\ar[d]\\
\bigoplus_{\sigma, \tau\in S_n} \bigotimes_{i=1}^{n}
\Ho(M^I_{f(i),h(\tau(\sigma(i)))})\ar[d]\\
\bigoplus_{\rho\in S_n} \bigotimes_{i=1}^{n}
\Ho(M^I_{f(i),h(\rho(i))})\ar[d]\\
\Co^{M}_{1,r}(f,h)}.\]
where the second arrow permutes the order in the tensor products, the third arrow is the product in
open string topology, and the other arrows are identities. The formula above shows that compositions
in $\Co^{M}_{1,r}$ are essentially determined by products in open string topology.\\

\begin{figure}[h]
\centering
\includegraphics[width=0.4\textwidth]{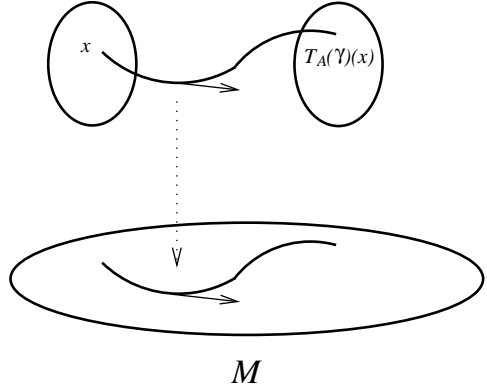} \caption{ \ Parallel transportation on a fiber bundle over $M.$}\label{876}
\end{figure}

Let $G$ be a compact Lie group and $\pi\colon P\longrightarrow M$
be a principal $G$-bundle over $M$. We let $\scr{A}_{P}$ be the
space of all  connections on $P$. There are many ways to think of
a connection on a principal fiber bundle, for us the most
important fact is that associated to such a connection  $A
\in \scr{A}_{P}$ there is a notion of parallel transportation, i.e.  if $\gamma \colon I \longrightarrow M$ then $A$
gives rise in a canonical way to a map
$$T_{A}(\gamma): P_{x(0)} \longrightarrow P_{x(1)}.$$
Figure
\ref{876} illustrates the process of parallel transportation.
The most important properties of the operators $T_{A}(\gamma)$ is
that it is independent of reparametrizations of the curve
$\gamma$, it depends continuously on both $A$ and $\gamma$, and if
$\gamma_1
\circ \gamma_2$ is the path obtained by the concatenation of path
$\gamma_1$ and $\gamma_2$ then
$$T_{A}(\gamma_1 \circ \gamma_2)=T_A(\gamma_2) \circ T_A(\gamma_1)
.$$

Our next goal is to prove the following result.

\begin{thm}\label{sl}
{\em There is a natural map $\Ho:\scr{A}_{P} \longrightarrow
\HL_{1,r}(M).$}
\end{thm}

For each connection $A \in \scr{A}_{P}$ we construct a functor
$$\Ho_{A}\colon\Co_{1,r}^{M}\longrightarrow \ve.$$ It sends an object $f$ of $\Co_{1,r}^{M}$
into
$$\Ho_{A}(f)=\Ho({P_{f}})=H({P_{f}})[\dim(\overline{f})],$$
where $P_{f(i)}$ denotes the restriction of $P$ to $f(i)\subseteq
M$ and $${P_{f}}=\prod_{i\in[n]} P_{f(i)}.$$ Theorem
\ref{sl} follows from the next result.

\begin{prop}\label{t6}
{\em The map $\Ho_{A}\colon \Co_{1,r}^{M}\longrightarrow
\ve$ sending $f$ into $\Ho_{A}(f)$ defines a one dimensional
restricted homological quantum field theory.}
\end{prop}
We need to define linear maps  $$\Ho_{A}\colon
\Co_{1,r}^{M}(f,g)
\longrightarrow  \ho(\Ho({P_{f}}),
\Ho({P_{g}})).$$  By the previous discussion an element of
$\Co_{1,r}^{M}(f,g)$ is a tuple $(\sigma,t)=(\sigma,t_1,...,t_n)$
where $\sigma \in S_n$ and $$t_i \in
\Ho(M^I_{f(i),g(\sigma(i))}).$$ The map
\[\Ho_{A}(\sigma, t)\colon \Ho(P_f)\to
\Ho(P_g) \] is defined as follows.
Consider the projection map $\pi\colon P_f \longrightarrow
\overline{f},$ and let $x$ be a chain $x\colon
K_{x}\longrightarrow P_f$ where $x=(x_{1},\dots ,x_{n}).$ The
domain of $\Ho_{\Lambda}(\alpha, t)(x)$ is given by
\[K_{\Ho_{A}(\alpha, t)(x)}=K_{x}\times_{\overline{f}}
\prod _{i\in[n]}K_{t_{i}}\] The map
$\Ho_{A}(\alpha,t)(x)\colon
K_{\Ho_{A}(\alpha,t)(x)}\longrightarrow P_{g}$ is given by
\[[\Ho_{\Lambda}(\alpha, t)(x)](y;s_{1},\cdots
s_{n})_i=[T_A(t_{i}(s_{i}))](x_{i}(y))\] where $y\in K_{x},
\ s_{i}\in K_{t_{i}}.$\\

\begin{figure}[h]
\centering
\includegraphics[width=0.5\textwidth]{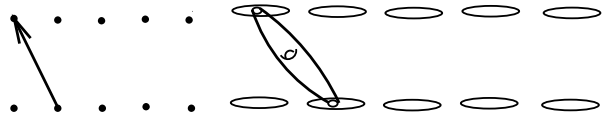}
\caption{ \ Motivation for homological matrices. }\label{matana5}
\end{figure}

The construction above produces objects of $\HL_{1,r}$ from
connections in principal bundles. It would be interesting to
determine what is the image under this map of known families of
connections, say for example flat connections or the  $N$-flat
connections introduced in \cite{angel}. We like to mention that
there is a remarkable analogy between \HL
\ in dimension one and the algebra of matrices. Recall \cite{dpgi}
that we can identify the space of matrices with the vector space
generated by bipartite graph with a unique edge with starting
point in $[n]$ and endpoint in $[m]$. For example Figure
\ref{matana5} represents on the left a $5\times 5$  matrix and on
the right a higher dimensional homological analogue.

\begin{figure}[h]
\centering
\includegraphics[angle=90,width=0.3\textwidth]{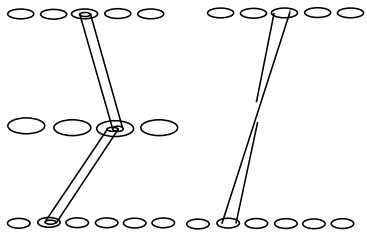}
\caption{ \ Composition for homological graphs. }\label{matana2}
\end{figure}

Thus it is clear how to define the higher dimensional homological
analogue of  the algebra of matrices;   we call the new algebra
the algebra of homological graphs. It is simple given by
\[\mathrm{HG}(f,g)=\des\bigoplus_{j\in[m], i\in[n]}\Ho(M^I_{f(j),g(i)})\]
where $f:[m] \longrightarrow D(M)$ and $g:[n] \longrightarrow
D(M).$ Moreover one can define a product on the space of
homological matrices that generalizes the usual product of
matrices. It is given by combining the usual matrix product with
the product of open strings, see \cite{cd3} for details. The
higher dimensional product is represented in Figure \ref{matana2}.

\begin{figure}[h]
\centering
\includegraphics[width=0.5\textwidth]{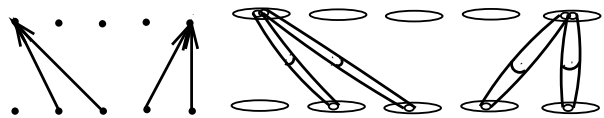}
\caption{ \ Higher dimensional morphisms in the $Schur_4$ category. }\label{matana4}
\end{figure}

\begin{figure}[h]
\centering
\includegraphics[width=0.4\textwidth]{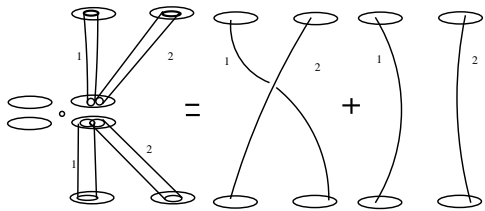} \caption{ \ Example of composition in the
$\mathrm{Schur}_2$ category. }\label{7777}
\end{figure}

Once we have defined an homological analogue of the algebra of
matrices, the problem of extending the usual constructions with
matrices to the higher dimensional case arises naturally. In
\cite{cd3} we explored that question  and found that several well-known
constructions for matrices may indeed be generalized to the
homological context. One of them is the possibility of  defining
homological Schur algebras and Schur categories. Recall that the
$Schur_k$ category
\cite{DP2} is such that its objects are positive integers and its
morphisms are given by
$$ \mathrm{Schur}_k(n,m)=Sym^k( End(\mathbb{C}^n, \mathbb{C}^m)),$$ i.e.  $Schur_k$ is the
$k$-symmetric power of the category of linear maps between the
vector spaces  $\mathbb{C}^n$. An example of a morphism in the
category $Schur_4$ is displayed on the left of Figure
\ref{matana4}. On the right there is an example of a morphisms in
the higher dimensional Schur category. The product rule in the
symmetric powers of algebras or categories where introduced in
\cite{DP2} and has been further studied in \cite{DP0, RDEP3, DR}.
Figure
\ref{7777} shows, schematically, an example of composition in the
$
\mathrm{Schur}_2$ category. Notice that in this case the product
of basis elements is not an element of the basis. Representations
of homological $\mathrm{Schur}_k(n,n)$ algebras are deeply related
with one dimensional homological quantum field theories
\cite{cd3}.

\section{Two dimensional homological quantum field theory }\label{ttt}

In this section we study homological quantum field theories in
dimension $2$. Our first goal is to generalize the map from
connections to $\HL_{1,r}$ to the $2$-dimensional situation. Our
second goal is two define a the membrane homology $\mathcal{H}(M)$
associated with each compact oriented manifold $M$. The matrix
graded algebra $\mathcal{H}(M)$ may be regarded as a
$2$-dimensional analogue of the Chas-Sullivan string
topology.\\

Let $M$ be a compact oriented smooth manifold and consider the
space $M^{S^1}$ of free loops on $M$. We assume that we are given
a complex Hermitian line bundle $L$ on  $M^{S^1}$. According to
Segal \cite{Se2} a $B$-field or  string connection on $L$ is a
rule that assigns to each pair $(\Sigma, y)$  where $\Sigma$ is a
surface with a boundary and $y$ is a map $y\colon
\Sigma
\longrightarrow M$ a parallel transportation operator
$$B_{y} \colon L_{\partial(\Sigma)_{-}} \longrightarrow
L_{\partial(\Sigma)_{+}},$$ where the extension of $L$ to
$(M^{S^1})^{n}$ is defined by the rule $$L_{(x_1,...,x_n)}=L_{x_1}
\otimes ... \otimes L_{x_{n}}.$$ The assignment $y \longrightarrow
B_y$ is assumed to have the following properties:

\begin{itemize}
\item{It is a continuous map taking values in unitary operators. Therefore we have induced maps
$B_{y}\colon L^{1}_{\partial(\Sigma)_{-}} \longrightarrow
L^{1}_{\partial(\Sigma)_{+}}$ between the corresponding circle
bundles.}

\item{It is transitive with respect to the gluing of surfaces.}

\item{It is a parametrization invariant.}
\end{itemize}

Let $\mathcal{B}_L$ be the space of $B$-fields or string
connections on $L$. Our next goal is to prove the following
result.
\begin{thm}
{\em There is a natural map  $\mathcal{B}_L
\longrightarrow  \HL_{2,r}(M).$}
\end{thm}

Thus for each $B$ field we need to construct a functor
$$\Ho_{B}\colon
\Co^{M}_{2,r} \longrightarrow \ve.$$ Since the only compact manifold without boundary of dimension $1$ is a circle, then
an object in $\Co^{M}_{2,r}$ is  a map $f
\colon [n] \longrightarrow D(M)$.  The functor $H_B$ is defined by the rule $$\Ho_{B}(f)=
\Ho({L}_{f}^{1})= H({L}_{f}^{1})[\dim(\overline{f})],$$ where $$\Ho({L}_{f}^{1})=
\Ho(L^{1}_{f(1)\times ... \times f(n)}).$$ The notation
$L^{1}_{f(1)\times ... \times f(n)}$ makes sense since
$$f(1)\times ... \times f(n) \subseteq M
\times ... \times M \subseteq M^{S^1} \times ... \times M^{S^1}.$$

\begin{prop}
{\em The map $\Ho_{B}\colon \Co_{2}^{M}\longrightarrow \ve$
sending $f$ into $\Ho(L_{f}^1)$ defines a two dimensional
restricted homological quantum field theory.}
\end{prop}

We need to define linear maps  $$\Ho_{B}\colon
\Co_{2,r}^{M}(f,g)
\longrightarrow  \ho(\Ho({L_{f}^1}),
\Ho({L_{g}^1})).$$
Suppose $f
\colon [n] \longrightarrow D(M)$, $g:[m] \longrightarrow D(M)$, that we are given  a chain
$x:K_x \longrightarrow L_f^1 $, and that we have another chain
$y:K_y \longrightarrow M_{f,g}^{\Sigma}$, where $\Sigma$ is a
surface with $n$ incoming boundaries and $m$ outgoing boundaries.
The maps $\pi(x)$ and $e_0(y)$ allow us to define the domain of
$\Ho_{B}(x)(y)$ as follows
$$K_{\Ho_{B}(x)(y)}=K_x \times_{\overline{f}}K_y. $$
The map
$$\Ho_{B}(x)(y): K_{\Ho_{B}(x)(y)}\longrightarrow L_g^1$$
is given by
$$\Ho_{B}(y)(x)(s,t)=B_{y(t)}[x(s)] .$$\\

Next we proceed to construct the membrane topology associated with
each compact oriented manifold $M$. Let us take a closer look at
objects in the category $\Co^{M}_{2,r}$. We focus our attention on
objects $f \colon [n] \longrightarrow D(M)$ such that $f$ is
constantly equal to $M$, thus the map $f$ becomes irrelevant and
this type of objects are indexed by positive integers. A morphism
from $[n]$ to $[m]$ is a homology class of the space $M^{\Sigma}$
of maps from $\Sigma$ into $M$ that are constant around the
boundaries of $\Sigma$, where $\Sigma$ is a compact oriented
surface with $n$ incoming boundary components and $m$ outgoing
boundary components. We shall further restrict our attention to
connected surfaces $\Sigma$. \\

\begin{figure}
\begin{center}
\includegraphics[width=0.4\textwidth]{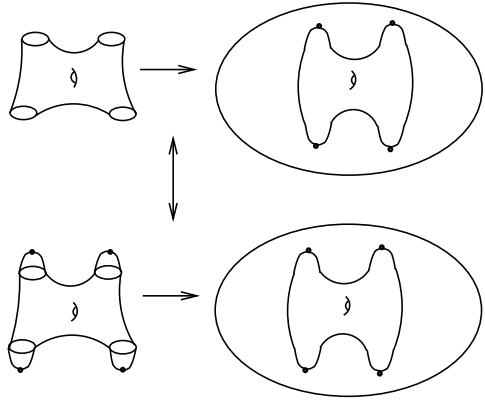}
\caption{ \ An element in $M^{\Sigma}$ and the corresponding element in $M^{\Sigma_{n}^{m}}$ . \label{bra1}}
\end{center}
\end{figure}

For integers $n,m \geq 1,$ let $\Sigma_{n,g}^{m}$ be a connected
Riemann surface of genus $g$ with $n$ incoming marked points and
$m$ outgoing marked points. Let $M^{\Sigma_{n,g}^{m}}$ be the
space of smooth maps from $\Sigma_{n,g}^{m}$ to $M$ which are
constant in a neighborhood of each marked point. If $\Sigma$ is a
genus $g$ surface with $n$ incoming boundary components and $m$
outgoing boundary components, then the spaces $M^{\Sigma}$ and
$M^{\Sigma_{n,g}^{m}}$ are homotopically equivalent, see Figure
\ref{bra1}  for an example illustrating the homotopy equivalence between $M^{\Sigma}$ and
$M^{\Sigma_{n,g}^{m}}$.  Therefore we have that
$$H(M^{\Sigma})=H(M^{\Sigma_{n,g}^{m}}).$$
We are going to use the following algebraic definition. We say
that an algebra $(A,p)$ is a matrix graded if $$A=
\bigoplus_{n,m=1}^{\infty} A_{n}^{m}, \ \ \ \ p \colon A_{n}^{m}\otimes A_{m}^{k}\to
A_{n}^{k}, \mbox{ \ \ \ and \ \ \ } p\mid_{A_{n}^{m}\otimes
A_{l}^{k}}=0 \mbox{\ \ \ if  \ \ }l\neq m. $$

\begin{defi} {\em The membrane homology of a compact oriented manifold $M$ is given by
$$\mathcal{H}(M)=\bigoplus_{n,m=1}^{\infty}
\Ho_{n}^{m}(M),$$ where $\Ho_{n}^{m}(M)=\bigoplus_{g=0}^{\infty} \Ho^{m}_{n,g}(M)$ and
$\Ho^{m}_{n,g}(M)=H(M^{\Sigma_{n,g}^{m}})[m\dim(M)].$}
\end{defi}

\begin{figure}[h]
\centering
\includegraphics[width=0.4\textwidth]{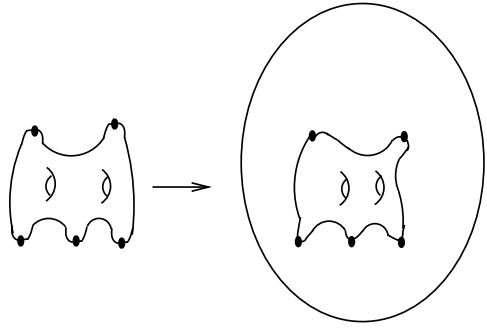} \caption{ \ Element in $M^{\Sigma_{3,2}^{2}}$.}\label{cobra2}
\end{figure}

\begin{figure}[h]
\centering
\includegraphics[width=0.4\textwidth]{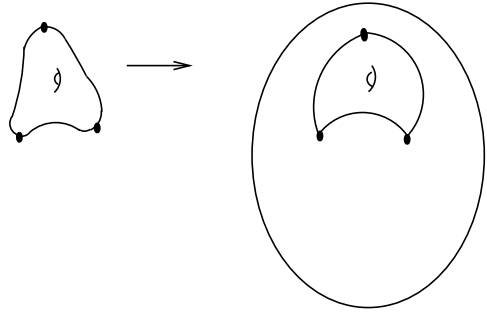} \caption{ \ Element in $M^{\Sigma_{2,1}^{1}}$.}\label{cobra1}
\end{figure}

Figure \ref{cobra2}  and Figure \ref{cobra1} show, schematically,
examples of an element in $\Ho^{2}_{3,2}(M)$ and an element of
$M^{\Sigma_{2,1}^{1}}$.

Using the composition rule in $\Co^{M}_{2,r}$ and the fact the
gluing of $\Sigma_{n,g}^{m}$ and $\Sigma_{k,g}^{l}$ is
$\Sigma_{n,g+k}^{l}$ one arrives to the following conclusion.

\begin{thm}
{\em $\h(M)$ is a matrix graded algebra.}
\end{thm}

\begin{figure}[h]
\centering
\includegraphics[width=0.4\textwidth]{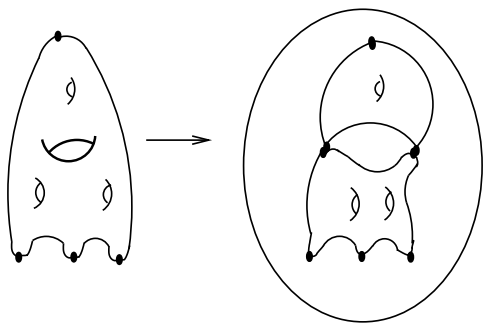}
\caption{ \ Element $i(x,y) \in M^{\Sigma_{3,4}^{1}}$. \label{cobra3}}
\end{figure}

The product of the elements shown in Figures \ref{cobra2}  and
\ref{cobra1} is shown in Figure \ref{cobra3}.\\

Finally  we like to mention that membrane homology of a manifold
$M$ comes equipped with a canonical representation. For a vector
space $V$ we let $T_{+}(V)=
\bigoplus_{n=1}^{\infty}V^{\otimes n}.$

\begin{thm}
{\em $T_{+}(\Ho(M))$ where $\Ho(M)=H(M)[dim M]$ is a
representation of $\h(M)$.}
\end{thm}

We have seen that the membrane topology is an interesting
algebraic structure associated to each oriented manifold. It would
be interesting to compute it explicitly for familiar spaces, and
also to study its relation with other types of two dimensional
field theories, such as topological conformal field theories in
the sense of \cite{kim1, kim2}.

\section{Conclusion}

In this work we introduced three new topological invariants for
compact oriented manifolds. The first invariant is a functor
$$Oman \longrightarrow tfd\mbox{-}alg$$ from $Oman$ the groupoid of compact
oriented manifolds into $tfd\mbox{-}alg$ the category of
transversal algebras over the operad $C(fD^d)$ of chains of framed
little $d$-discs. The functor is given by the correspondence
$$M \longrightarrow \Ca(M^{S^{d}}),$$
which maps a compact oriented manifold  $M$ into the space
$\Ca(M^{S^{d}})$ of chains of maps from the $d$-sphere into $M$.
The second invariant depends on the choice of a compact oriented
manifold $Y$. It is a functor $$Oman \longrightarrow
t1\mbox{-}Cat$$ from compact oriented manifolds into
$t1\mbox{-}Cat$ the category of small transversal $1$-categories,
i.e. categories over the operad $\Ca(I)$ of chains of little
intervals; the functor is given by the correspondence
$$M \longrightarrow \Ca(M^{S(Y)}),$$
sends a compact oriented manifold into the $1$-category
$\Ca(M^{S(Y)}).$ Our third invariant depends only on the choice of
an integer $d\geq 1$.  It is a functor $$Oman \longrightarrow
g\mbox{-}Cat$$ from $Oman$ to $g\mbox{-}Cat$  the category of
small graded categories, and its given  by the correspondence
$$M \longrightarrow \HL_{d,r}(M),$$ sending a manifold $M$ into the
category  $\HL_{d,r}(M)$ of restricted homological quantum fields
theories on $M$. Given a compact oriented manifold we have
constructed several examples of objects in the category
$\HL_{d,r}(M)$. We paid especial attention to the cases $d=1$ and
$d=2$. In the former case we see that connections on line bundles
on $M$ are sources of objects in $\HL_{1,r}(M)$. Likewise
$B$-fields on line bundle over $M^{S^1}$ is a source of examples
of object in $\HL_{2,r}(M)$. From the notion of homological
quantum field in dimension $2$ we constructed an algebra
associated to each compact oriented manifold called the membrane
topology of $M$. This algebra may be
thought as a $2$-dimensional generalization of the string topology of Chas and Sullivan.\\

Finally, let us mention a few open problems and ideas for future
research that arise naturally from the results of this work:

\begin{itemize}
\item Further examples of HLQFT are  needed. A potential source of
examples could be the higher dimensional generalizations of
$B$-fields, for example using the higher-dimensional notion of
parallel transport of Gomi and Terashima \cite{GT}, or perhaps the
parallel transport for $n$-Lie algebras recently developed in
\cite{SU}.

\item The main obstacle towards an explicit description of the
category of homological quantum fields theories is that only for a
handful of spaces the homology groups $H(M^L)$ are known
explicitly. Results along this line are very much welcome.

\item It would be interesting to investigate to what extend the
notion of HLQFT can be extended to yield topological invariants
for singular (non-smooth) manifolds. A step forward in that
direction have been taken by Lupercio, Uribe and Xicotencatl in
\cite{Ur3} where they consider string topology on orbifolds.
\end{itemize}

\noindent  edmundocastillo@gmail.com\\
\noindent Escuela de Matem\'aticas, Facultad de Ciencias, Universidad Central de Venezuela,
Av. Los Ilustres Los Chaguaramos, A.P.: 20513, Caracas 1020, Venezuela.\\

\noindent ragadiaz@gmail.com\\
\noindent Grupo de F\'isica-Matem\'atica, Universidad Experimental de las Fuerzas Armadas, Caracas 1010, Venezuela.\\
\end{document}